\documentclass[11pt]{amsart}%%%%%,leqno
\usepackage{amsmath,amssymb, graphics, amscd,latexsym }

\makeatletter

\newtheorem{Theorem}{Theorem}%[section]
\newtheorem{Lemma}[Theorem]{Lemma}
\newtheorem{Corollary}[Theorem]{Corollary}
\newtheorem{Proposition}[Theorem]{Proposition}

\theoremstyle{remark}
\newtheorem{Remark}[Theorem]{Remark}

\newtheorem{Example}[Theorem]{Example}

\begin{document}
\newcommand{\eps}{\varepsilon}
\newcommand{\om}{\omega}
\newcommand\Om{\Omega}
\newcommand\la{\lambda}
\newcommand\vphi{\varphi}
\newcommand\vrho{\varrho}
\newcommand\al{\alpha}
\newcommand\La{\Lambda}
\newcommand\si{\sigma}
\newcommand\be{\beta}
\newcommand\Si{\Sigma}
\newcommand\ga{\gamma}
\newcommand\Ga{\Gamma}
\newcommand\de{\delta}
\newcommand\De{\Delta}

\newcommand\cA{\mathcal  A}
\newcommand\cB{\mathcal B}
\newcommand\cD{\mathcal  D}
\newcommand\cM{\mathcal  M}
\newcommand\cN{\mathcal  N}
\newcommand\cT{\mathcal  T}
\newcommand\cP{\mathcal  P}
\newcommand\cp{\mathcal  p}
\newcommand\cQ{\mathcal  Q}
\newcommand\cG{\mathcal G}
\newcommand\cq{\mathcal  q}
\newcommand\cc{\mathcal  c}
\newcommand\cs{\mathcal  s}
\newcommand\cS{\mathcal  S}
\newcommand\ct{\mathcal  t}
\newcommand\cZ{\mathcal  Z}
\newcommand\cR{\mathcal  R}
\newcommand\cu{\mathcal  u}
\newcommand\cU{\mathcal  U}
\newcommand\co{\mathcal  o}
\newcommand\cO{\mathcal  O}
\newcommand\cv{\mathcal  v}
\newcommand\cV{\mathcal  V}
\newcommand\cx{\mathcal  x}
\newcommand\cX{\mathcal  X}
\newcommand\cw{\mathcal  w}
\newcommand\ck{\mathcal  k}
\newcommand\cK{\mathcal  K}
\newcommand\cW{\mathcal  W}
\newcommand\cz{\mathcal  z}
\newcommand\cy{\mathcal  y}
\newcommand\ca{\mathcal  a}
\newcommand\ch{\mathcal  h}
\newcommand\cH{\mathcal  H}
\newcommand\cF{\mathcal F}

\newcommand\bfC{\mbox {\bf  C}}
\newcommand\bfN{\mbox {\bf  N}}
\newcommand\bfT{\mbox {\bf  T}}
\newcommand\bfP{\mbox {\bf  P}}
\newcommand\bfp{\mbox {\bf  p}}
\newcommand\bfQ{\mbox {\bf  Q}}
\newcommand\bfq{\mbox {\bf  q}}
\newcommand\bfc{\mbox {\bf  c}}
\newcommand\bfs{\mbox {\bf  s}}
\newcommand\bfS{\mbox {\bf  S}}
\newcommand\bft{\mbox {\bf  t}}
\newcommand\bfZ{\mbox {\bf  Z}}
\newcommand\bfR{\mbox {\bf  R}}
\newcommand\bfu{\mbox {\bf  u}}
\newcommand\bfU{\mbox {\bf  U}}
\newcommand\bfo{\mbox {\bf  o}}
\newcommand\bfO{\mbox {\bf  O}}
\newcommand\bfv{\mbox {\bf  v}}
\newcommand\bfV{\mbox {\bf  V}}
\newcommand\bfx{\mbox {\bf  x}}
\newcommand\bfX{\mbox {\bf  X}}
\newcommand\bfw{\mbox {\bf  w}}
\newcommand\bfk{\mbox {\bf  k}}
\newcommand\bfK{\mbox {\bf  K}}
\newcommand\bfW{\mbox {\bf  W}}
\newcommand\bfz{\mbox {\bf  z}}
\newcommand\bfy{\mbox {\bf  y}}
\newcommand\bfa{\mbox {\bf  a}}
\newcommand\bfh{\mbox {\bf  h}}
\newcommand\bfH{\mbox {\bf  H}}

\newcommand\bbC{\mbox {\mathbb C}}
\newcommand\bbN{\mbox {\mathbb N}}
\newcommand\bbT{\mbox {\mathbb T}}
\newcommand\bbP{\mbox {\mathbb P}}
\newcommand\bbQ{\mbox {\mathbb Q}}
\newcommand\bbS{\mbox {\mathbb S}}
\newcommand\bbZ{\mbox {\mathbb Z}}
\newcommand\bbR{\mbox {\mathbb R}}
\newcommand\bbU{\mbox {\mathbb U}}
\newcommand\bbO{\mbox {\mathbb O}}
\newcommand\bbV{\mbox {\mathbb V}}
\newcommand\bbX{\mbox {\mathbb X}}
\newcommand\bbK{\mbox {\mathbb K}}
\newcommand\bbW{\mbox {\mathbb W}}
\newcommand\bbH{\mbox {\mathbb H}}

\newcommand\apeq{\fallingdotseq}
\newcommand\Lrarrow{\Leftrightarrow}
\newcommand\bij{\leftrightarrow}
\newcommand\Rarrow{\Rightarrow}
\newcommand\Larrow{\Leftarrow}
\newcommand\nin{\noindent}
\newcommand\ninpar{\par \noindent}
\newcommand\nlind{\nl \indent}
\newcommand\nl{\newline}
\newcommand\what{\widehat}
\newcommand\tl{\tilde}
\newcommand\wtl{\widetilde}
\newcommand\order{\mbox{\text{order}\/}}
\newcommand\GL{\text{GL}\/}
\newcommand\Spec{\text{Spec}\/}
\newcommand\weight{\text{weight}\/}
\newcommand\ord{\text{ord}\/}
\newcommand\Int{\text{Int}\/}
\newcommand\grad{\text{grad}\/}
\newcommand\Ind{\text{Ind}\/}
\newcommand\Disc{\text{Disc}\/}
\newcommand\Ker{\text{Ker}\/}
\newcommand\Image{\text{Image}\/}
\newcommand\Coker{\text{Coker}\/}
\newcommand\Id{\text{Id}\/}
\newcommand\id{\text{id}}
\newcommand\dsum{\text{\amalg}}
\newcommand\val{\text{val}}
\newcommand\minimum{\text{minimum}\/}
\newcommand\modulo{\text{modulo}\/}
\newcommand\Aut{\text{Aut}\/}
\newcommand\PSL{\text{PSL}\/}
\newcommand\Res{\text{Res}\/}
\newcommand\rank{\text{rank}\/}
\newcommand\codim{\text{codim}\/}
\newcommand\Cone{\text{Cone}\/}
\newcommand\maximum{\text{maximum}\/}
\newcommand\Vol{\text{Vol}\/}
\newcommand\lcm{\text{lcm}\/}
\newcommand\degree{\text{degree}\/}
\newcommand\Pol{\cal {POL}}
\newcommand\Ts{Tschirnhausen}
\newcommand\TS{Tschirnhausen approximate}
\newcommand\Stab{\text{Stab}\/}
\newcommand\civ{complete intersection variety}
\newcommand\nipar{\par \noindent}
\newcommand\wsim{\overset{w}{\sim}}
\newcommand\Gr{\bfZ_2*\bfZ_3}
\newcommand\QED{~~Q.E.D.}
\newcommand\bsq{$\blacksquare$}
\newcommand\bff{\mbox {\bf  f}}
\newcommand\newcommandby{:=}
\newcommand\inv{^{-1}}
\newcommand\nnt{(\text{nn-terms})}
\renewcommand{\subjclassname}{\textup{2000} Mathematics Subject Classification}
\def\mapright#1{\smash{\mathop{\longrightarrow}\limits^{{#1}}}}
%\def\maprightt#1#2{\smash{\mathop{\longrightarrow}\limits^{#1}}}
%左矢印
\def\mapleft#1{\smash{\mathop{\longleftarrow}\limits^{{#1}}}}
%同下向き矢印
\def\mapdown#1{\Big\downarrow\rlap{$\vcenter{\hbox{$#1$}}$}}
\def\mapdownn#1#2{\llap{$\vcenter{\hbox{$#1$}}$}\Big\downarrow\rlap{$\vcenter{\hbox{$#2$}}$}}
%同上向き矢印
\def\mapup#1{\Big\uparrow\rlap{$\vcenter{\hbox{$#1$}}$}}
\def\mapupp#1#2{\llap{$\vcenter{\hbox{$#1$}}$}\Big\uparrow\rlap{$\vcenter{\hbox{$#2$}}$}}

%同右下向き
\def\rdown#1{\searrow\rlap{$\vcenter{\hbox{$#1$}}$}}
\def\semap#1{\searrow\rlap{$\vcenter{\hbox{$#1$}}$}}

%同右上向き
\def\rup#1{\nearrow\rlap{$\vcenter{\hbox{$#1$}}$}}
\def\nemap#1{\nearrow\rlap{$\vcenter{\hbox{$#1$}}$}}
%同左下向き
\def\ldown#1{\swarrow\rlap{$\vcenter{\hbox{$#1$}}$}}
\def\swmap#1{\swarrow\rlap{$\vcenter{\hbox{$#1$}}$}}
%同左上向き
\def\lup#1{\nwarrow\rlap{$\vcenter{\hbox{$#1$}}$}}
\def\nwmap#1{\nwarrow\rlap{$\vcenter{\hbox{$#1$}}$}}
\def\defby{:=}
\def\eqby#1{\overset {#1}\to =}
\def\inv{^{-1}}
\def\bnu{{(\nu)}}

\title[ Fundamental group of  sextics of torus type
%{\it Draft: \today}
]
{  Fundamental group of sextics of torus type  }

\author[
%Normally smooth divisors  {\it Draft: \today}
M. Oka {\tiny and} D.T. Pho]
{Mutsuo Oka {\tiny and} Duc Tai Pho}
\address{\vtop{
%\hbox{Mutsuo Oka}
\hbox{Department of Mathematics}
\hbox{Tokyo Metropolitan University}
\hbox{1-1 Mimami-Ohsawa, Hachioji-shi}
\hbox{Tokyo 192-0397}
\hbox{\rm{E-mail}: \vtop{\hbox{\rm oka@comp.metro-u.ac.jp}\hbox{\rm pdtai@math.metro-u.ac.jp}}}}}

%\today
%\thanks{}

\begin{abstract}
We show that the fundamental group of the complement of any irreducible
tame torus sextics in 
$\bfP^2$ is isomorphic to $\bfZ_2*\bfZ_3$ except one class.
The exceptional class has the configuration of the singularities $\{C_{3,9},3A_2\}$
and the fundamental group is bigger than $\bfZ_2*\bfZ_3$.
In fact, the Alexander polynomial is given by 
$(t^2-t+1)^2$.
For the proof, we first reduce the assertion to maximal curves and 
then we  compute the fundamental groups for maximal tame torus curves.
% We also show that the assertion does not hold for non-tame
% torus curves, by a simple example.
\end{abstract}
\maketitle

\pagestyle{headings}
\section{Introduction}
 Recall that a sextic curve
$C: F(X,Y,Z)=0$ in $\bfP^2$
is called a {\em torus curve of type
(2,3)} if there exists an expression: $F(X,Y,Z)=
F_2(X,Y,Z)^3+F_3(X,Y,Z)^2$ where $F_2(X,Y,Z)$, $F_3(X,Y,Z)$ are 
polynomials of degree 2 and 3 respectively.
Let $C_2$ and $C_3$ be the conic and the cubic
defined by $F_2(X,Y,Z)=0$ and $F_3(X,Y,Z)=0$ respectively.
$C$ is called a {\em tame} torus curve if the singularities are only on
the intersection $C_2\cap C_3$.

In \cite{Pho}, the second author classified every possible singularities
on tame torus sextics of type (2,3). In particular, he showed 
that there exist 7 moduli of the maximal tame torus curves and
their configurations of the singularities are given as follows.
\begin{eqnarray}
\{C_{3,15}\},\{C_{9,9}\},\{B_{3,10},A_2\}, \{Sp_1,A_2\},
\{B_{3,8},E_6\}, \{C_{3,7},A_8\}, \{C_{3,9},3A_2\}
\label{max-list}\end{eqnarray}
For a configuration $\Si$ in (\ref{max-list}) we denote the corresponding moduli space
of sextics of torus type by
$\cM(\Si)$. It is also shown   that each of the
moduli space is connected. Thus the topology of the complement
$\bfP^2-C$
is independent of the choice of a generic curve $C\in \cM(\Si)$.
Any irreducible tame torus curve can be degenerated into one of them and 
 maximal tame torus curves are rational curves.

The purpose of this paper is to study of the geometry 
of these curves from the viewpoint of the fundamental group
of the complement of curves.

In \S 2, we prepare a key lemma (Lemma \ref{isomorphism lemma}) which reduces the
computation
of the fundamental groups to the case of maximal curves.

In \S 3, we  compute the dual curves of the maximal curves.
 Some of the maximal curves have a big singularity like
$C_{3,7}, C_{3,9}, C_{3,15}$ which are not locally irreducible but have 
the same tangent cone. We show the self-duality
of these singularities.
% using a method to compute the intersection 
%numbers of the dual images  of the respective components.

In \S 4, we give
 the main result of this paper (Theorem \ref{main-theorem}) which state:
\[\pi_1(\bfP^2-C)\cong \bfZ_2*\bfZ_3\]
for any irreducible tame torus curve 
and the Alexander polynomial is given by $t^2-t+1$ with 
 one  exceptional class $\cM(\{C_{3,9},3A_2\})$.
For a curve $C$ in this exceptional moduli space,
the fundamental group is represented by
\begin{eqnarray}
 \pi_1(\bfC^2-C)&\cong & \langle g_1,g_2,g_3 | \{g_i,g_j\}=e, i\ne j, (g_3g_2g_1)^2=
(g_2g_1g_3)^2=(g_1g_3g_2)^2\rangle\\
 \pi_1(\bfP^2-C)&\cong & \langle g_1,g_2,g_3 | \{g_i,g_j\}=e, i\ne j, (g_3g_2g_1)^2=1\rangle
\end{eqnarray}
and the Alexander polynomial is given by $(t^2-t+1)^2$.

%There are more cases where 
%the isomorphism $\pi_1(\bfP^2-C)\cong \bfZ_2*\bfZ_3$ does not hold for
%non-tame sextics.
In \S 5, we give an  example of a non-tame torus curve $C$ with three
 $E_6$
 and an $A_1$ 
singularities so that 
the fundamental group of this curve is isomorphic to 
$B_4(\bfP^1)$
 (Theorem \ref{non-tame}).

\section{Zariski pencil and Puiseux parametrization}
\subsection{Zariski-pencil}
Let $C$ be a curve of degree $d$ in $\bfP^2$. To compute the fundamental
group $\pi_1(\bfP^2-C)$,
it is usually most convenient to use the Zariski's pencil method, 
which we  recall briefly.
First choose a  line $L_\infty$
and choose a base point $b$ on $ L_\infty$.
The family of lines $L_\eta$ which pass through $b$ is called 
{\em the pencil at $b$}.
Unless otherwise stated, we take $b$ as the base point of the
fundamental group as well. A pencil line $L_\eta$ is called 
{\em singular} if $L_\eta\cap C$ contains less than $d$ points.
%This is the case if $L_\eta$ passes through a singular point of $C$ 
%or it is tangent to a simple point of $C$.
Take   homogeneous coordinates $(X,Y,Z)$ so that 
$L_\infty$ is defined by $Z=0$ and $b=(0,1,0)$ for simplicity.
Let $W=\{\eta_1,\dots,\eta_\ell\}$ be the set of the parameters
corresponding to the  singular pencil lines by the correspondence
$\eta\leftrightarrow L_\eta\maketitle=\{X-\eta Z=0\}$.
Choose a generic pencil line $L_{\eta_0}$. Thus $\eta_0\in \bfC-W$.
 By the local triviality of the line section $L_\eta\cap C$
over $\bfC-W$ , 
$\pi_1(\bfC-W,\eta_0)$ acts on $\pi_1(L_{\eta_0}-L_{\eta_0}\cap C)$.
This gives the monodromy relations
\[\{ g= g^\si; g\in \pi_1(L_{\eta_0}-L_{\eta_0}\cap C),\si\in
\pi_1(\bfC-W,\eta_0)\}\]
 The  van Kampen-Zariski theorem says that
$\pi_1(\bfP^2-C)$ 
%(resp. $\pi_(\bfC^2-C)$) are
is  isomorphic to the 
quotient group of $\pi_1(L_{\eta_0}-L_{\eta_0}\cap C)$ 
%(resp. $\pi_1((L_{\eta_0}\cap \bfC^2-L_{\eta_0}\cap C)  $ 
 by the 
monodromy relations.

The affine fundamental group $\pi_1(\bfC^2-C)$  can be computed by the
exact same way, replacing $\pi_1(L_{\eta_0}-L_{\eta_0}\cap C)$
by 
$\pi_1(L_{\eta_0}^a-L_{\eta_0}^a\cap C)$ where 
$L_{\eta_0}^a:=L_{\eta_0}-L_{\eta_0}\cap L_{\infty}$. Moreover 
for a generic line at infinity, we have a
central exact sequence ({\cite{Central}):
\[
1\to \bfZ\to \pi_1(\bfC^2-C)\mapright{\iota} \pi_1(\bfP^2-C)\to 1
\]
The generator $1\in \bfZ$ is represented by a lasso of $L_\infty$ and it is
homotopic to the ``big circle'' $\ell\circ S_R\circ \ell\inv$
where $S_R$ is the big circle  $|y|=R,~R\gg 1$
 in $L_{\eta_0}^a-L_{\eta_0}\cap C$ which contains all intersection points
$C\cap L_{\eta_0}$ and $\ell$ is a path joining $S_R$ and the base point.
\begin{Remark}
%In general, there are exactly $d(d-1)$ singular pencil lines counting 
%the multiplicity. 
In practice,
it is extremely difficult to read the monodromy relations
for curves which are defined over $\bfC$. Suppose
that we are interested in the fundamental $\pi_1(\bfP^2-C)$,
where $C$ has a prescribed configuration of the singularities $\Si$
and $\degree~ C=n$ is fixed. Let $\cM(\Si;n)$ 
be the moduli space of curves of degree $n$ with configuration
of the singularities $\Si$.
As the fundamental group does 
not change if we move the curve in a connected component of  the moduli space
$\cM(\Si)$, we are free to move the 
original curve in this component. So if possible,
it is convenient to choose a curve defined over the real numbers $\bfR$
which has as many singular pencil lines over $\bfR$ as possible.
%This is yet not sufficient.
 Usually a choice of a suitable pencil also  makes the computation 
  easier. \end{Remark}
\subsection{Puiseux parametrization}
Let $h(t)=\sum_{i=0}^\infty a_i t^i$  be a convergent power series
with complex coefficients. Let $n$ be a given positive integer.
The {\em first characteristic power} $P_1(h(t);n)$ is defined by
the integer
$\min\{j>0;a_j\ne 0,j\not \equiv 0~\mod ~n\}$
(see \cite{Spain}). Let $\nu_1:=P_1(h(t);n)$ and $n^{(1)}:=\gcd(n,\nu_1)$.
If $\nu_1<\infty$, we define $\nu_2:=P_1(h(t);n^{(1)})$ and
$n^{(2)}:=\gcd(n^{(1)},\nu_2)$ and so on. As the integers $n,n^{(1)},\dots$ are
decreasing, they become stationary after a finite steps. So we assume that 
$n^{(k-1)}>n^{(k)}=n^{(k+1)}$. We put
\[\cP(h(t);n):=\{\nu_1,\dots,\nu_k\},\quad
D(h(t);n):=\{n^{(1)},\dots, n^{(k)}\}\]
%Putting $n^{(0)}=n, \nu_{k+1}=\infty$,
%$h_j(t):=\sum_{\nu_j\le i<\nu_{j+1}} a_it^i$ is in fact a power series of 
%$t^{n_{j}}$ and $\val_t h_j(t)=\nu_j$ for $j=0,\dots, k$.

Let $(C,O)$ be a germ of an irreducible curve
with Puiseux pairs $\{(m_1,n_1),\dots,(m_k,n_k)\}$.
Recall that $\gcd(m_i,n_i)=1$ and $m_i>m_{i-1}n_i$ for 
$i=1,\dots,k$ with $m_0=1$.
Let $(x,y)$ be coordinates so that  $y=0$ defines the tangent cone.
Then $C$ can be parametrized as
$x(t)=t^N$ and 
$y=\phi(t)=\sum_{i=S}^\infty a_it^i$ so that 
\begin{eqnarray}
&\cP(\phi(t);N)=\{m_1n_2\cdots n_k,m_2n_3\cdots n_k,\dots, m_k\},\\
&D(\phi(t);N)=\{n_2\cdots n_k,n_3\cdots n_k,\dots, n_k\}
\end{eqnarray}
where $N=n_1\cdots n_k$.
Put $S=\val_t \phi(t)$. Note that $S\le m_1n_2\cdots n_k$
and~ $S\equiv 0\mod N $ if and only if 
$S< m_1n_2\cdots n_k$. The number 
$s:=S/N$ is called {\em the Puiseux order}
 of $y(x^{1/N})$ in \cite{Oka-sextics}.
Recall that $2\le s\le m_1/n_1$ and $s$ is an integer  if
$s<m_1/n_1$.
%The original function $f(x,y)$ is recovered (up to a multiplication of a unit)
%as 
%\[
%f(x,y)=\prod_{j=0}^{N-1}(y-\phi(t \omega^j)),\quad \omega:=\exp(2\pi i/N)\]
Consider an irreducible curve $C$ which is described as above.
To see the behavior of the intersection $C\cap \{y=\eta\}$, we wish to express
$x$ as a function of
$y$.
This is the case when we compute the fundamental group $\pi_1(\bfP^2-C)$
using the pencil $\{y=\eta;\eta\in \bfC\}$.  For
this purpose, we take the new parameter
$\tau$ which is defined by
$\phi(t)=\tau^{S}$.  Then we have
\begin{Lemma}\label{inverse-Puiseux}
We can write  $C$ as $y=\tau^S$ and $x=\psi(\tau)=\tau^N\psi_0(\tau)$ so that
$\psi_0(0)\ne 0$ and 
\begin{eqnarray*}
&\cP(\psi_0(\tau);N)&= 
\{ m_1n_2\cdots n_k-S,m_2n_3\cdots n_k-S,\dots, m_k-S\},\quad S<m_1n_2\cdots n_k\\
&\cP(\psi_0(\tau);\frac N{n_1})&=\{m_2n_3\cdots n_k-S,\cdots,
m_k-S\} ,\quad S=m_1n_2\cdots n_k
\end{eqnarray*}
\end{Lemma}
\begin{proof}
First, we can write $\phi(t)=t^S\phi_0(t)$ with $\phi_0(0)\ne 0$.
Then 
\[
\begin{cases} &\cP(\phi_0(t);N)=\{m_1n_2\cdots n_k-S,\dots, m_k-S\},\quad
S<m_1n_2\cdots n_k\\
            &\cP(\phi_0(t);N/n_1)=\{m_2n_3\cdots n_k-S,\dots,m_k-S\}, 
\quad S=m_1n_2\cdots n_k
\end{cases}
\]
Thus by \cite{Spain}, Lemma 5.1,  and Lemma 5.2,
the assertion follows immediately.
\end{proof}
Suppose that we have to compute the (local or global) fundamental group of
a curve which  have a reducible singularity at the origin $O$.
The above lemma plays an important role to compute the local or global
fundamental group in such a situation. See \S 4.

\begin{Example}\label{local-strategy}
Assume that $s=2$ and $C$ has  a single Puiseux pair $\{(m,n)\}$ with 
$m/n>2$ and $y=0$ to be the
tangent cone. Then 
$C$ is parametrized as 
\begin{eqnarray*}\label{standard}
x=t^n,\quad y=a_{2}t^{2n}+\cdots+a_{\frac mn} t^m+\cdots,\quad a_2,a_{\frac
mn}\ne 0
\end{eqnarray*}
and putting $y=\tau^{2n}$, we get another   parametrization
\begin{eqnarray*}\label{reverse}
y=\tau^{2n},\quad
x=\tau^n(b_0+b_{1}\tau^n+\cdots+b_{k}\tau^{kn}+b_{\frac mn-2} \tau^{m-2n}+\cdots),
\quad k=[m/n]-2
\end{eqnarray*}
and $b_0,b_{\frac mn-2}\ne 0$.
The local topological behavior depends only on $\{s,m,n\}$.
We remark  that for the line $x=\eta,|\eta|\ll 1$ intersects with $C$ at  $n$
points locally near $O$, while the  line $y=\eta,|\eta|\ll 1$ 
intersects $C$ at $2n$ points.
\end{Example}
\subsection{Isomorphism theorem}
%We prepare an isomorphism  theorem  which we use later.
Let $G$ be a group. The commutator subgroup is denoted by $D(G)$.
The first homology of $G$ is by definition the quotient group
$G/D(G)$ and we denote it by $H_1(G)$. A free group of rank $n$ is 
denoted by $F(n)$
\begin{Lemma}\label{isomorphism lemma}
Let $G$ be a group such that  $D(G)$ is a free group $F(n)$ 
with  $n <\infty$. Suppose that  we have 
a surjective endomorphism $\vphi: G\to G$ which induces 
an isomorphism   $\bar \vphi$ on $H_1(G)$.
Then $\vphi$ is an isomorphism.
\end{Lemma}
\begin{proof}
This is  observed in \cite{Oka-sextics}. First we consider the
surjective  homomorphism
$\vphi':=\vphi|_{D(G)}:D(G)\to D(G)$. By the assumption,
we have an isomorphism
$D(G)/\Ker(\vphi')\cong D(G)$. If  $\Ker(\vphi')$ is not
trivial,
we get an obvious contradiction  to the Hopfian property of the free groups
%\[n>\rank(D(G)/\Ker(\vphi'))=\rank(D(G))=n\]
 (see Theorem  2.13,
\cite{M-K-S}). Now the assertion follows from  Five Lemma:
\[
\begin{matrix}
1&\to&D(G)&\to&G&\to &H_1(G)&\to &1\\
&&    \mapdown{\vphi'}&&     \mapdown{\vphi}&&\mapdown{\bar\vphi}&&\\
1&\to&D(G)         &\to&    G       &\to& H_1(G)&\to &1
\end{matrix}
\] 
\end{proof}
As is well-known, $D(\bfZ_p*\bfZ_q)$ is a free group of rank $(p-1)(q-1)$,
 we obtain 
\begin{Corollary}\label{isomorphism}
Assume that $\vphi:\bfZ_p*\bfZ_q\to \bfZ_p*\bfZ_q$ is a surjective
homomorphism which gives an isomorphism on the first homology.
%where $p,q$ are  relatively coprime integers . 
Then $\vphi$ is an
isomorphism.
\end{Corollary}
\subsection{Geometric homomorphism}Suppose that $C,C'$ are reduced  curves of the same degree
and assume that the line at infinity $L_\infty:=\{Z=0\}$ is generic for 
$C$ and $C'$. A homomorphism 
$\phi: \pi_1(\bfC^2-C)\to \pi_1(\bfC^2-C')$ is {\em geometric}
if it preserves 
``the big circles''  so that it induces an homomorphism 
$\bar\phi:\pi_1(\bfP^2-C)\to \pi_1(\bfP^2-C')$ and  the commutative
diagrams:
\begin{eqnarray}\label{geometric diag}
\begin{matrix}
1&\to  & \bfZ&\to &\pi_1(\bfC^2-C)&\mapright{\iota}&\pi_1(\bfP^2-C)&\to&
1\\
&&\mapdown{\Id}&&\mapdown{\phi}&&\mapdown{\bar\phi}&&\\
1&\to  & \bfZ&\to &\pi_1(\bfC^2-C')&\mapright{\iota'}&\pi_1(\bfP^2-C')&\to&
1
\end{matrix}
\end{eqnarray}
%We call a homomorphism $\psi:\pi_1(\bfP^2-C)\to \pi_1(\bfP^2-C)$
%{\em geometric} if it is induced by a geometric homomorphism
%$\phi: \pi_1(\bfC^2-C)\to \pi_1(\bfC^2-C')$.
where $\iota,\iota'$ are canonical homomorphisms induced by the
respective inclusion maps.
By the definition, we have 
\begin{Proposition} Assume that $\phi:\pi_1(\bfC^2-C)\to \pi_1(\bfC^2-C')$
is a geometric homomorphism. In  (\ref{geometric diag}), $\phi$
is an isomorphism if and only if $\bar\phi$ is an isomorphism.
\end{Proposition}

Assume that we have a degeneration family of plane curves  $C_t\, (|t|\le 1)$
so that $C_t$ is reduced for any $t$ and $C_t, t\ne 0$, 
has the same configuration of singularities. In such a situation we have a
canonical surjective homomorphism
$\phi:\pi_1(\bfC^2-C_0)\to \pi_1(\bfC^2-C_t)$
which is geometric.
\begin{Corollary}\label{bigger}
Let $C$ be a torus sextic curve. If $\pi_1(\bfC^2-C)$ is not isomorphic
to the braid group $B(3)$, then  $\pi_1(\bfP^2-C)\not \cong \bfZ_2*\bfZ_3$. 
\end{Corollary}
\section{Dual of maximal sextics}
The information for the dual curves often plays an crucial 
role for the study of singularities on a plane curve of a given degree.
For example, we have used the dual curve information to show the
impossibility of the 
degeneration $\{C_{3,9},3A_2\}\to \{B_{3,12},A_2\}$ in \cite{Pho}.
In this section, we study the dual curves of maximal sextics.
\subsection{Dual singularities}
We consider the set of germs of  irreducible curves
$(C,O)$ which has  given Puiseux pairs
 $\cP:=\{(m_1,n_1),\dots, (m_k,n_k)\}$ and $y=0$ is the tangent line of $(C,O)$.
This set is denoted by 
$\si(\cP)$.
Let $s$ be the Puiseux order.
The subset of $\si(\cP)$ with the Puiseux order $s$
 is denoted by $\si(\cP;s)$ 
in \cite{Oka-sextics}. Then 
$\si(\cP)$ is the union of $\si(\cP;j)$ for $j=2,\dots, [m_1/n_1]$ and 
$m_1/n_1$. The family
$\{\si(\cP;s); s=2,\dots, [m_1/n_1],\,\text{and }\,s=m_1/n_1\}$ is  called
{\em the flex stratification} of $\si(\cP)$. Take
$(C,O)\in \si(\cP;s)$. Then $(C,O)$ is parametrized as
\[
 x=\xi(t),\quad y=\phi(t),\quad\xi(t):= t^N,~
 \phi(t):=\sum_{i=sN}^\infty a_i
t^i,\quad N :=n_1\cdots n_k
\]
%Recall  that $s>1$ and   $y=0$ is the
%tangent cone of $(C,O)$.
Though the parameter is local,
the coordinates $(x,y)$ is the affine coordinates
$x=X/Z,y=Y/Z$, unless otherwise stated.
  Recall that the defining function of 
$(C,O)$ is  given 
\[f(x,y)=\prod_{j=1}^N (y-\phi(t \omega^j)),
\quad
\omega=\exp(\frac{2\pi i}{N})\]
where we  replace $t^N$ by $x$.
Using the dual homogeneous coordinates $(U,V,W)$ of $(X,Y,Z)$
and the affine coordinates $u:=U/V$ and $w:=W/V$, we can parametrize
the dual curves $(C^*,O^*)$  as  (see
\cite{Oka-sextics})
\[
u(t):=-\frac{\phi'(t)}{\xi '(t)}=-\sum_{i=sN}\frac{ia_i}{N} t^{i-N},\quad
w(t):=\frac{\xi (t)\phi'(t)-\xi '(t) \phi(t)}{\xi '(t)}=\sum_{i=sN}
\frac{(i-N)a_i}{N} t^{i}
\]
Note that $\text{val}_tu(t)=(s-1)N$.
In the homogeneous coordinates, $O^*$ corresponds to $(0,1,0)$.
Taking the parameter $\tau$ so that 
$u(t(\tau))=\tau^{(s-1)N}$, we can write 
$w\maketitle=\psi(\tau)$ and   $\psi(\tau)\colon=\sum_{j} b_j \tau^j$ and 
the characteristic powers $\cP(\psi(\tau);(s-1)N)$
is completely described by $s$ and the Puiseux pairs $\cP$ (\cite{Oka-sextics}).
We want to describe the dual singularities of certain reducible 
singularities. For this purpose, we recall  the argument 
% in \cite{Oka-sextics}
as follows.
%For the general case, we have the following. 
Put $\ell=[m_1/n_1]$
and write 
\[\phi(t):= a_2t^{2N}+\cdots+a_{\ell}t^{\ell N}+
a_{\frac{m_1}{n_1}}t^{m_1n_2\cdots n_k}+\dots\]
Put 
\begin{eqnarray*}
&t(\tau)&=\tau (c_0+c_1\tau^N+\cdots+c_{\ell-2}\tau^{(\ell-2)N}+
c_{\frac{m_1}{n_1}-2}t^{(m_1-2)n_2\cdots n_k}+\cdots)\\
&\psi(\tau)&=b_2\tau^{2N}+\cdots+b_{\ell} \tau^{\ell N}+b_{\frac{m_1}{n_1}}\tau^{m_1N/n_1}
+\cdots
\end{eqnarray*}
\begin{Proposition}\label{generic-dual-Puiseux}
 Suppose that $s=2$.
Then   $c_0^{N}=-\frac{1}{2a_{2}}$ and 
the following coefficients $c_0,c_1\dots$ are inductively determined 
by the equality $u(t(\tau))=\tau^{(s-1)N}$.
The  coefficients $c_1,\dots,c_p$ of $t(\tau)$  and thus
the coefficients  $b_{2},\dots, b_{2+p}$ of
$\psi(\tau)$  depend only on
the first  $p+1$ coefficients
$a_{2},\dots, a_{2+p}$  of $\phi(t)$ and the choice of $c_0$
for  $p\le m_1/n_1-2$. The coefficient
$b_{2}$ is given by $\frac 1{4a_2}$.
\end{Proposition}

 For our later purpose, we recall
two special cases. See \cite{Oka-sextics} for detail.

\noindent
(1)  ({\bf Self-duality of the generic stratum}) Suppose that $s=2$
and $(C,O)\in \si(\cP;2)$.
Then we have $(C^*,O^*)\in \si(\cP;2)$. Two special cases
which we use later:

\noindent
(a) $\cP=\emptyset$. In this case, $(C,O)$ is a smooth germ and 
$N=1$ and the parametrization takes the form:
$y=\sum_{i=2}^{\infty} a_i x^i$, $ a_2\ne 0$.
The dual curve is parametrized as 
$w=\sum_{i=2}^{\infty} b_i u^i, b_2\ne 0$
%  with $b_2=\frac 1{4a_2}$ 
and the coefficients
$b_p$ depend only on $a_2,\dots, a_p$. Note that $C^*$ is a regular (not a flex)
point of $C^*$.

\noindent
(b) $\cP=\{(m,n)\}$ with $m/n>2$.
The parametrization   takes the following form:
\[x=t^n,\quad y=\sum_{i=2}^{\ell} a_{i} t^{in}+a_{\frac mn} t^m+\cdots,
\quad \ell:=[m/n]\]
with $a_{2},a_{\frac mn}\ne 0$. We can easily see that $f(x,y)$ takes
the following form:
\begin{eqnarray*}
f(x,y)=y_1^n-(a_{\frac mn})^n x^m+\text{higher terms},\quad
y_1=y-(a_{2} x^2+\cdots a_{\ell} x^{\ell})
\end{eqnarray*}
This implies that $f$ is non-degenerate with respect to 
$(x,y_1)$ and $(C,O)=B_{m,n}$. See \S 3.2  for the definition.
The dual curve $(C^*,O^*)$ has a similar parametrization
\begin{eqnarray}\label{first-coeff}
x=t^n,\quad y=\sum_{i=2}^{\ell} b_{i} t^{in}+b_{\frac mn} t^m+\cdots,
\quad
b_{2n}=\frac 1{4a_2}
\end{eqnarray}

\noindent
(2) Non-generic case with $s=m_1/n_1$. In this case,
$C$ is parametrized as 
\[ x=t^N,\quad y=\phi(t), ~~\phi(t)=a_{m_1n_2\cdots n_k}t^{m_1N/n_1}+\cdots\]
and 
the dual singularity  $(C^*,O^*)$ is given by
\begin{eqnarray*}
u = \tau^{(m_1-n_1)N/n_1},\quad
w=\psi(\tau),\quad \psi(\tau)= b_{m_1n_2\cdots n_k}\tau^{m_1N/n_1}+\cdots
\end{eqnarray*}
and the Puiseux order is $\frac{m_1}{m_1-n_1}$.
As for the Puiseux pairs,  we have
\begin{eqnarray}\cP(\psi(\tau))=\begin{cases}
&\{ (m_1,m_1-n_1),(m_2,n_2),\dots,(m_k,n_k)\},\quad m_1-n_1>1\\
&\{(m_2,n_2),\dots,(m_k,n_k)\}, \qquad m_1-n_1=1
\end{cases}\end{eqnarray}
\subsection{Dual singularity of a reducible singularity}
Recall that we have introduced the following topological equivalent
classes of curve singularities which appear on irreducible tame torus curves
(\cite{Pho}):
\begin{eqnarray}
\begin{cases}\label{further-types}
B_{p,q}:  \: y^p+x^q=0  \: (\text{Brieskorn-Pham type})\\
C_{p,q}:  \: y^p+x^q+x^2y^2=0,\quad
\frac 2p+\frac 2q\le 1\\
%D_{4,7}:  \: y^4+x^3y^2+x^7=0\\
Sp_1:  \: (y^2-x^3)^2+(xy)^3=0\\
%Sp_2:  \: (y^2-x^3)^2-y^6=0
\end{cases}
\end{eqnarray}
Now we  consider a reducible curve with a common tangent cone.
For example, $B_{p,q}, p<q$ has $\gcd(p,q)$ irreducible components with $x=0$ as the
tangent cone. $C_{3,q},\, q\ge 7$,
 has 2 (respectively 3) irreducible components with $y=0$ as
the tangent cone if $q$ is   odd (resp. even).

 Assume that two germs of plane curve singularities at the
origin
$(C_1,O)$ and $(C_2,O)$  have the same tangent cone. 
The topological equivalence class of 
$(C_1\cup C_2,O)$ is determined by
the respective Puiseux pairs and the intersection number $I(C_1,C_2;O)$.
We assume for simplicity that they have at most one Puiseux pairs
$\{(m,n)\},~m>2n$ and $\{(m',n')\},~m'>2n'$ respectively
and their Puiseux orders are 2.
By abuse of notation, we understand that $C_1$ is smooth if $n=1$.
For our purpose, we only need to consider  the case $n=n'$ or $n=1$.
Let $f_1(x,y),
f_2(x,y)$ be the respective defining functions and suppose that  their Puiseux
parametrizations be given as
\begin{eqnarray}
&C_1: x=t^{n},&\quad y=\phi_1(t),\quad
\phi_1(t)=a_{2}t^{2n}+\dots+a_{\frac mn}t^m+\cdots\\
&C_2:x=t^{n'},&\quad y=\phi_2(t),\quad
\phi_2(t)=a_{2}'t^{2n'}+\dots+a_{\frac{m'}{n'}}'t^{m'}+\cdots
\end{eqnarray}
Now  the intersection number is given as
\[
I(C_1,C_2;O)=\val_t f_2(t^n,\phi_1(t))\]
Assume that $\alpha\in \bfQ$ be the minimum
of $j$ such that 
$a_{j}\ne a_{j}'$. We consider the case 
$\alpha\le \min\{m/n,m'/n'\}$.
As $f_2(x,y)$ is written as
\begin{eqnarray*}
f_2(x,y)=y_2^{n'}-(a_{\frac{m'}{n'}})^{n'} x^{m'}+\text{(higher terms)},\quad
y_2= y-(a_{2n'} x^2+\cdots+a_{[\frac{m'}{n'}]} x^{n'[m'/n']}),
\end{eqnarray*}
we can easily compute that
\begin{eqnarray}
I(C_1,C_2;O)=\alpha n n'
\end{eqnarray}
Under the same assumption, we consider the parametrization of the dual curves:
\begin{eqnarray}
&C_1^*: x=t^{n},& \quad y=\psi_1(t),\quad
\psi_1(t)=b_{2}t^{2n}+\dots+b_{\frac mn}t^m+\cdots,\quad b_2=\frac{1}{4a_2}\\
&C_2^*:x=t^{n'},&\quad y=\psi_2(t),\quad
\psi_2(t)=b_{2}'t^{2n}+\dots+b_{\frac{m'}{n'}}'t^m+\cdots,\quad b_2'=\frac{1}{4a_2'}
\end{eqnarray}
Assume that $n=n'$. By Proposition \ref{generic-dual-Puiseux}
we have $b_{j}=b_{j}'$ for $j<\alpha$ and $b_{\alpha }\ne b_{\alpha }'$.
Thus we get $I(C_1^*,C_2^*;O^*)=I(C_1,C_2;O)=\alpha nn'$.

Assume that $n=1$ and $a_2\ne a_2'$. Then by a similar argument,
 $I(C_1,C_2;O)=2n'$
and $I(C_1^*,C_2^*;O^*)=2n'$.
\subsection{Dual of singularities on tame torus curve (local).}
First we consider reducible germs with a common tangent cone.
It is shown that $C_{3,n},~n=7,8,9,12,15$,  appears as a singularity
of irreducible  tame torus curves (\cite{Pho}).
  A common
nature of these singularities is the following. 

\noindent
$(\sharp)$:  There are a smooth component $L$
corresponding to the face supporting 
$y^3+x^2y^2$.  The face supporting $x^2y^2+x^n$ 
gives  an $A_{n-3}$ singularity  $(K,O)$.  
$K$  is  irreducible (respectively $K$ has two smooth components
 $K_1$ 
and $K_2$)
 if $n$ is odd  (resp.  $n$ is even).
In any case, the tangent cone is $y=0$. The intersection
number
is given by $I(L,K;O)=4$ if $n$ is odd and 
$I(L,K_i;O)=2$ for $i=1,2$ and $I(K_1,K_2;O)=n/2-1$  for $n$ even.

Conversely $C_{3,n}$ is characterized by this
property. For example, assume that $n=2m$ and $(C,O)$ has three smooth
components at $O$, which satisfy the above intersection criterion.
Take an analytic coordinates $(x_1,y_1)$ so that $K_1=\{y_1=0\}$.
Then $K_2$ is defined by an analytic function of the form 
$j_2(x_1,y_1)=y_1+a_{m-1} x_1^{m-1}+\dots$ with $a_{m-1}\ne 0$ and $L$ is
defined  by an analytic function of the form  $j_0(x,y)=y+b_2x_1^2+\dots$
with
$b_2\ne 0$. Thus $(C,O)$ is defined, using  a new coordinates
$(x_1,y)$ with $y_1=y+\eps x^{m-1}, ~\eps+a_{m-1}\ne 0$ by the function
\begin{eqnarray*}
j(x,y)=&j_0(x,y+\eps x^{m-1})(y+\eps x^{m-1})j_2(x,y+\eps x^{m-1})\\
=&y^3+b_2 x^2y^2+b_2\eps( a_{m-1}+\eps) x^n+\text{(higher terms)}
\end{eqnarray*}

Now suppose that we have a tame torus curve $C$ with a singularity at $O$.
We consider the dual singularity $(C^*,O^*)$.
In the local classification argument, we have shown that 
 the Puiseux order of each  component $L,K$ or $L,K_1,K_2$
is  2 with respect to the  fixed affine
coordinates $(x,y)$.

 For example, 
 let us consider  a tame torus curve $C$ with $(C,O)\cong C_{3,15}$
and we assume that $y=0$ is the
common tangent cone. We have shown in \cite{Pho}
that $C$ is defined by a polynomial $f(x,y)$ which is   written as
$f(x,y)=f_2(x,y)^3+f_3(x,y)^2$  where the conic $f_2(x,y)=0$ is tangent 
with $y=0$ at $O$ and the cubic $f_3(x,y)=0$
has a node at $O$ and one branch is tangent to $y=0$
and $I(f_2,f_3;O)=6$.
 The   component
$K$ has a
$A_{12}$ singularity
and we have the following parametrization.
\begin{eqnarray}
L: y& =& a_2 x^2+a_3x^3+\cdots\\
K: x &=& t^2,\quad y=\sum_{i=2}^6 a_{i}' t^{2i}+a_{13/2} t^{13}+\cdots
\end{eqnarray}
and $a_2\ne a_2'$. Thus $\alpha=2$ and   $I(L,K;O)=4$.
Thus $(L^*,O^*)$ is smooth and 
$(K^*,O^*)$ is again a generic $(2,13)$ cusp
and $I(L^*,K^*;O^*)=4$.  This implies that 
$(C^*,O^*)\cong C_{3,15}$.  In an exact same discussion, we see that 
$(C_{3,n},O)$ are self dual for any $n\ge 7$.

Other singularities on irreducible tame torus curves with a
common tangent and having $s=2$ as the Puiseux order are 
$A_{3n-1} ~(n=1,\dots,6)$ and $B_{3,n}~(n=6,8,10,12)$.
By the same discussion, we can see that these   singularities
%on  tame torus curves
 are self dual. Thus in conclusion, we have
\begin{Proposition}The following singularities on tame 
torus curve
are self-dual. 
\[ C_{3,n},n=7,\dots,15; \quad A_{3n-1},n=1,\dots, 6;
\quad B_{3,2m},m=3,\dots, 6\]
\end{Proposition}
Secondly we consider reducible germs with several tangent cones.
 Singularities with several tangent cones  which we have in mind are
$C_{6,6}, C_{6,9}, C_{6,12}, C_{9,9}$. Each of
these singularities has two components in the  tangent cones. Let $B_x$ 
(respectively 
$B_y$) be the  union of irreducible components  which has $x=0$ (resp.  $y=0$)
 as   the  tangent cone.
 Each of the irreducible
component is generic and therefore on  the dual curve and  each of them is
isomorphic to the original one. Let $O_x:=(1,0,0)$ and $O_y=(0,1,0)$.
Then $(C_{m,n},O)^*=(B_x^*,O_x)\cup (B_y^*,O_y)$ and 
they are self dual for $m,n\ge 6$: $(B_x^*,O_x)\cong (B_x,O)$ and
$(B_y^*,O_y)\cong (B_y,O)$.  For example, consider 
$C_{6,9}$.  Note that $(B_x,O_x)\cong A_3$ and $(B_y,O_y)\cong A_6$.
%A tame torus singularity $(C,O)$ is isomorphic to $C_{6,9}$  if the conic 
%$C_2$ consists of two lines passing
%at $O$ and 
% the cubic $C_3$ is   nodal   at $O$ for which one tangent is a line 
%of $C_2$. Assume that $C_2$ has  $y=0$ as a component and 
%the tangent cones of the node of
%$C_3$ at $O$ is defined by $xy=0$. Then $(C,O)\cong C_{6,9}$. 
On the dual curve
$C^*$, the dual singularity splits into a $A_3$ and $A_6$
with a common tangent line.

Finally we consider
other singularities on torus curves. Exceptional singularities are 
$A_2, E_6$, $B_{4,6}$ and $Sp_1$. As is well-known, the dual of a cusp
$B_{n,n+1}$ is a flex point of order $n-1$.  This implies that the dual
of
$B_{4,6}$ is 
$A_{5}$.   As two components have a flex point at $O^*$,
 their Puiseux orders are 3. $Sp_1$ is an irreducible singularity
 with two Puiseux pairs $\{(2,3), (2,9)\}$ and the Puiseux order is $3$.
The generic form of $Sp_1$ is given by (up to $\PSL(3,\bfC)$-action)
\[(axy)^3+(y^2-x^3+a_{21}x^2y+a_{12}xy^2+a_{03}y^3)^2=0\]
and Puiseux pairs are $\cP=\{(3,2),(9,2)\}$. It has the parametrization
\[x=t^4,\quad y=\phi(t),
\quad \phi(t)=t^6+b_8 t^8+b_9 t^9+\cdots,b_9\ne 0\]
Thus $(Sp_1^*,O^*)$ has the parametrization
\[u=\tau^2,\quad w=\psi(\tau),\quad
\psi(\tau)=b_6\tau^6+b_8 \tau^8+b_9\tau^9+\cdots\]
Namely  $(Sp_1^*,O^*)$ is an $A_8$ singularity which has
the Puiseux order  3 and thus does not belong to  the generic
stratum.
\subsection{Generic dual curves (global)}
First we recall that irreducible tame torus curves can be degenerated into 
a  rational curve with one of the following configurations
(\cite{Pho}).
\begin{eqnarray}\label{maximal-curves}
\{C_{3,15}\}, \{C_{9,9}\}, \{C_{3,7},A_8\},
\{B_{3,10},A_2\}, \{Sp_1,A_2\}, \{B_{3,8},E_6\},
 \{C_{3,9},3A_2\}
\end{eqnarray}
%We study the dual curves of these maximal curves.
Let $\Si$ be one of the above configuration and let $\cM(\Si)$
be the set of tame torus curves with configuration $\Si$.
We say  that $C\in \cM(\Si)$ is generic if $C^*$ 
has only $A_2$ and $A_1$ singularities besides the singularities
which are dual to those of $C$. For a generic curve $C\in \cM(\Si)$,
the number of cusps and nodes are constant, we have observed that 
$\cM(\Si)$ is connected for each of the above 
maximal configuration.
The following is immediate from  the previous consideration 
and by computing the dual curve of an explicit generic curve.

\begin{Proposition} Let $\Si$ be one of the configuration in
(\ref{maximal-curves}) and let $C\in \cM(\Si)$  be a generic  curve. Let $n^*$ 
be the degree of $C^*$ and let $\Si^*$ be the  configuration of the
singularities of $C^*$. Then the dual curves are described by the following
table.
\begin{table}[htb]
\begin{tabular}[t]{|c|c|c|c|}
\hline
No&$\Si$ & $n^*$ & $\Si^*$ \\ \hline
1&$\{C_{3,15}\}$       & $9$ & $\{C_{3,15}, 9A_2,9A_1\}$    \\
2&$\{C_{9,9}\}$        & $8$ & $\{C_{3,7},A_8, 6A_2,5A_1\}$ \\
3&$\{ C_{3,7},A_8\}$   & $8$ & $\{C_{3,7},A_8, 6A_2,5A_1\}$ \\
4&$\{B_{3,10},A_2\}$   & $6$ & $\{B_{3,10}, A_2\}$          \\
5&$\{Sp_1,A_2\}$       & $6$ & $\{A_8, 3A_2,3A_1\}$         \\
6&$\{B_{3,8},E_6\}$    & $6$ & $\{B_{3,8}, 2A_2, A_1\}$     \\
7&$\{C_{3,9},3A_2\}$   & $6$ & $\{C_{3,9},3A_2\}$           \\
\hline
\end{tabular}
\end{table}
\end{Proposition}
\begin{Remark}
It is interesting to observe that in the last four cases, the dual curves
are 
also sextics. It is possible to show that they are torus curves.
In the case 5 and  6, the dual torus curves are not tame as it has $A_1$
singularity.
%Further detail will be studied elsewhere.
\end{Remark}
\section{Fundamental group of tame torus curves}
The fundamental group of the complement of a generic sextics of torus
type 
is isomorphic to $\bfZ_2*\bfZ_3$ by Zariski \cite{Za1}. The same
assertion is true for a certain 
 class of non-generic sextics of torus curves (\cite{Oka-sextics}).
The main result of this paper is:
\begin{Theorem}\label{main-theorem}
Let $C$ be an irreducible tame torus
sextic of type (2,3) defined by  the homogeneous polynomial $F(X,Y,Z)$ 
and let $M$ be the Milnor fiber of  $F(X,Y,Z)$. 

{\rm (1)} If $C\notin \cM(\{C_{3,9},3A_2\})$,
then $\pi_1(\bfP^2-C)\cong\bfZ_2*\bfZ_3$ and the generic Alexander
 polynomial is $t^2-t+1$ and the first Betti number $b_1(M)$ of $M$ is 2.

{\rm(2)} (Exceptional moduli)
For $C\in \cM(\{C_{3,9},3A_2\})$, we have
\begin{eqnarray}
 \pi_1(\bfC^2-C)&\cong & \langle g_1,g_2,g_3 ~|~ \{g_i,g_j\}=e, i\ne j, (g_3g_2g_1)^2=
(g_2g_1g_3)^2=(g_1g_3g_2)^2\rangle\\
 \pi_1(\bfP^2-C)&\cong & \langle g_1,g_2,g_3 ~|~ \{g_i,g_j\}=e, i\ne j, (g_3g_2g_1)^2=e\rangle
\end{eqnarray}
and the generic Alexander polynomial is given by $(t^2-t+1)^2$
and the first Betti number $b_1(M)$ is 4.
This moduli space $\cM(\{C_{3,9},3A_2\})$ is self dual.
\end{Theorem}
Recall that the {\em generic Alexander polynomial} of a curve $C$ is
defined by the Alexander polynomial
of $\pi_1(\bfC^2-C)$ with respect to a generic line at infinity.
The Milnor fiber $M$ is defined by the 
affine surface $\{(X,Y,Z)\in \bfC^3; F(X,Y,Z)=1\}$
and it is  the total space of the cyclic covering of order 6,
$p\colon M\to \bfP^2-C $, which is defined by
the quotient map of the $\bfZ/6\bfZ$-action induced by the monodromy map
$h\colon M\to M$.
We have seen that there are 7 configurations of  the possible maximal tame torus
curves, as is listed in (\ref{maximal-curves}).
In \cite{Pho}, it is  observed  that each of the
moduli space of the maximal configuration is an Zariski-open subspace of an affine
space
and thus it is  connected. 
Thus the topology of the complement is independent of the choice of a curve.
%Furthermore  any tame torus curve can be degenerated into one of
%the maximal torus curve, listed in (\label{maximal-curves}). 
We have
seen in
\cite{Pho} that there are degenerations
of curves of tame torus curves corresponding to 
\[
 \{C_{3,7},3A_2\} \to \{C_{3,8},3A_2\}\to \{C_{3,9},3A_2\}
\]
Any non-maximal tame torus curve $C$ can be degenerated into a
maximal  curve with configuration in (\ref{maximal-curves})
and if $C$ degenerates into 
$\{C_{3,9},3A_2\}$, it can be factored with  a degeneration 
into 
$\{C_{3,8},3A_2\}$.

For the proof of the assertion (1) of  Theorem \ref{main-theorem},
we prepare the following.
\begin{Lemma} \label{reduction to maximal curves} Let $C$ be an irreducible  tame
torus curve and let
$C_t$ be a degenerating family such that $C_1=C$ and $C_t,~t\ne 0$ has the same 
configuration of singularities with that of $C_1$.
Assume that $\pi_1(\bfP^2-C_0)\cong \bfZ_2*\bfZ_3$. Then $\pi_1(\bfP^2-C)\cong
\bfZ_2*\bfZ_3$.
\end{Lemma}
\begin{proof}
In fact, the degeneration family gives a surjective homomorphism
$\alpha: \pi_1(\bfP^2-C_0)\to \pi_1(\bfP^2-C)$ which is an isomorphism on the
first homology. Now take another    degeneration families
  $D_t$ such that
 $ D_t,t\ne 0$ is a generic torus curve with 6 $A_2$
and $D_0=C$. Such a family always exists.
We  get a  surjective homomorphism $\beta:\pi_1(\bfP^2-C)\to \pi_1(\bfP^2-D_t)$
  which also induces an isomorphism on the first
homology group.
Note that $\pi_1(\bfP^2-D_t)\cong\bfZ_2*\bfZ_3$ by \cite{Za1}.
 Now we
apply Lemma
\ref{isomorphism lemma} to conclude that the composition 
\[\bfZ_2*\bfZ_3\cong
\pi_1(\bfP^2-C_0)\mapright{\alpha}\pi_1(\bfP^2-C)\mapright{\beta}
\pi_1(\bfP^2-D_t)\cong\bfZ_2*\bfZ_3
\]
is an isomorphism. This implies both of  $\alpha, \beta$ are isomorphisms.
\end{proof} Thus the proof of Theorem \ref{main-theorem} is reduced to 
 the assertion for the maximal curves and 
a curve $C\in \cM(\{C_{3,8},3A_2\})$.
The assertion for these curves will be proved by
direct computations choosing   good members from respective moduli spaces.
To show that $\pi_1(\bfP^2-C)\cong \bfZ_2*\bfZ_3$ for 
 $C\in \cM(\Si)$
where $\Si$ is in 
(\ref{maximal-curves}) and $\Si\ne\{C_{3,9},3A_2\} $
or $\Si=\{C_{3,8},3A_2\}$, we  use the following lemma.
\begin{Lemma}\label{max-equivalence} The following conditions are equivalent
for a  tame torus curve $C$ and a generic line at infinity.
\begin{enumerate}
\item $\pi_1(\bfP^2-C)\cong \bfZ_2*\bfZ_3$.
\item %With respect to a generic line at infinity,
$\pi_1(\bfC^2-C)\cong B(3)$.
\item
There is a surjective homomorphism
$\psi: \bfZ_2*\bfZ_3\to \pi_1(\bfP^2-C)$
which gives  an isomorphism on the first homology.
\item There exists a surjective homomorphism 
$\phi :   F(2) \to D(\pi_1(\bfP^2-C))$. 
\end{enumerate}
\end{Lemma}
\begin{proof}Using a degeneration family $D_t$ from the generic torus curves as 
in the proof of Lemma \ref{reduction to maximal curves},
we have always a surjective homomorphism $\beta:\pi_1(\bfP^2-C)\to
\bfZ_2*\bfZ_3$.
Thus the assertion follows  from Lemma \ref {isomorphism lemma}.
\end{proof}

By Lemma \ref{max-equivalence}, to show $\pi_1(\bfP^2-C)\cong \bfZ_2*\bfZ_3$
for a  given sextic curve $C$
of torus type, we only need to show the existence of 
 generators $g_1,g_2$ given by lassos which satisfy
the braid relation $\{g_1,g_2\}=e$ and the torsion relation 
$(g_1g_2)^3=e$ 
%(see the next subsection). 
Usually we need  a few monodromy relations to get
these relations. Once we get these relations,
 we can ignore the other  monodromy relations
and we do not need any further computation.

Hereafter the base point of the fundamental group is the base point of the pencil
lines which we use. The local singularities in the  three configurations
$\{B_{3,10},A_2\}, \{Sp_1,A_2\},\{B_{3,8},E_6\}$  are 
irreducible.
%  and the computation is relatively easy.
\subsection{Notations and choice of the pencil}
Let $C$ be an irreducible torus curve of degree $6$.
A {\em lasso} is a loop represented as 
$\ell\circ \si\circ \ell\inv$ where $\si$ is a  loop, given by 
the boundary of a small normal disk of a regular point of $C$
and $\ell$ is a path joining $\si$ and the base point.
Thus through the Hurewicz homomorphism
$\xi: \pi_1(\bfP^2-C)\to H_1(\bfP^2-C)\cong \bfZ/6\bfZ$,
a lasso is mapped to the canonical generator of the first homology.
Recall that $\bfZ_2*\bfZ_3$ has two representations:
\begin{eqnarray*}
\bfZ_2*\bfZ_3&=&\langle \rho,\xi ~|~ \{\rho\xi\}=e,\, (\rho\xi)^3=e\rangle\\
&=&\langle a,b | a^2=b^3=e\rangle
\end{eqnarray*}
where $e$ is the unit element and
$\{\rho,\xi\}:=\rho\xi\rho(\xi\rho\xi)\inv$. Thus it implies
the braid relation $\rho\xi\rho=\xi\rho\xi$.
In the first representation, $\rho,\xi$ can be represented by lassos.
In the following figures,  we denote, for simplicity of drawing pictures,
 a
small lasso oriented in the counter clockwise direction by a bullet with 
a  path as in \cite
{Two,Oka-sextics}. 
Thus 
${\vrule height3pt depth-2.3pt width1cm} \bullet$ indicates
%$ {\vrule height3pt depth-2.3pt width1cm} \odot$.
$ {\vrule height3pt depth-2.3pt width1cm}\hspace{2pt} \cdot \hspace{-8.5pt}\circlearrowleft$.
%if we take $y=t$ as pencil lines.
%

For the computation of $\pi_1(\bfP^2-C)$ 
of  a  maximal tame torus curve $C$, we note that 
it has a big singularity $\xi\in C$ and it is 
usually better (by experiment) to choose Zariski pencil 
so that the pencil line, say $L_\rho$, passing 
through $\xi$ is the tangent cone of the singularity.
In this way, the monodromy relation  around $L_\rho$ contains more
relations.
Unless otherwise stated, we use $y=t$ as 
pencil lines  in this section. 
We explain how to read the monodromy relation by
 Example \ref{local-strategy}  in \S 2 with one Puiseux pair
putting $(m,n)=(13,2)$  which appears as a component of $C_{3,15}$.
\begin{eqnarray}\label{local-topo} 
y=\tau^{4},\,
x=b_1\tau^2+b_{2}\tau^4+\cdots+b_{5}\tau^{10}
+b_{\frac {11}2} \tau^{11}+\cdots
\end{eqnarray}
The local topology of $C$  is determined by the two terms
$b_1\tau^2,b_{\frac {11}2} \tau^{11} $
and the other terms does not change the topology.
For brevity,  we denote  the other terms  by {\it nn-terms}.
Then the parametrization (\ref{local-topo}) is simply written as
\[
 y=\tau^{4},\quad
x=b_1\tau^2+b_{\frac {11}2} \tau^{11}+\nnt
\]  
Let us consider a generic pencil $L_\varepsilon, 0<\varepsilon\ll 1$.
Assume $g_1,g_2,\dots,g_6$ are generators of the fundamental
group $\pi_1 (L_{\eta_0}-C)$.
% where  $\eta_0 =\eps$, $d$ is degree of the curve $C$.
 We may assume that $g_1,g_2,\dots,g_{4}$ are the
generators which correspond to those points bifurcated from $O$ of $C$. 
When $y=\eps\exp(\theta i)$ moves around the origin once, 
each branch of $\tau=\eps^{1/4}\exp(\theta i/4)$
moves  an arc of angle $\pi/2$.  The  topological
behavior of $4$ points among
 $C\cap \{y=\eps\exp(\theta i)\}$ looks like the movement of the two planets which accompany 
two respective satellites.
Two planets (which correspond to the term $b_1\tau^2$), 
do the half turn around the sun 
(=the origin). 
For a fixed $\sqrt y=\tau^2$, there exists $2$ roots of $\tau=\sqrt{\tau^{2}}$
which correspond to $2$ satellites.
%such that $\tau^n=\sqrt y$.
They do $\frac{11}{4}$-turns around
the respective planet.
This interpretation is useful when the local singularity is reducible.

Given a polynomial $f(x,y)$ we will denote by
$\Delta_x(f)$  
the discriminant of $f$ as a polynomial in $x$. Then $\Delta_x(f)$ is a 
polynomial in $y$, we denote $P_y$ is the set of roots of  
 $\Delta_x(f)=0$ in $y$. Thus $\lambda\in P_y$ 
corresponds to the  singular pencil line $y=\lambda$.

\subsection{Irreducible singularities} We first consider 3  moduli spaces,
$\{B_{3,10},A_2\}, \{Sp_1,A_2\}$, $\{B_{3,8},E_6\}$,
 which contain only irreducible singularities.
Hereafter $\eps$ is assumed to be  a sufficiently small
generic positive number throughout the paper.

\vspace{.2cm}  
\noindent
(I) {\bf Moduli space} $\cM(\{B_{3,10},A_2\})$. 
Let us consider the following curve $C_I\in \cM(\{B_{3,10},A_2\})$ which is defined by
\begin{eqnarray*}
&&C_I:f= (y(1-y)-x^2)^3+\left(y^2(1-y)-x^2y+xy^2+\frac{18}{25}y^3\right)^2\\
&&\Delta_x(f)=cy^{23}(949y-625)^3h_4(y), \quad P_y=\{0,625/949,
\text{roots of}~ h_4(y)=0\}
\end{eqnarray*}
where $h_4(y)=135771071y^4-676592025y^3+1149291875y^2-820856875y+284765625$.

\begin{figure}[htb]
\setlength{\unitlength}{1bp}
\begin{picture}(255,130)(-100,0)
\put(-150,130){\includegraphics{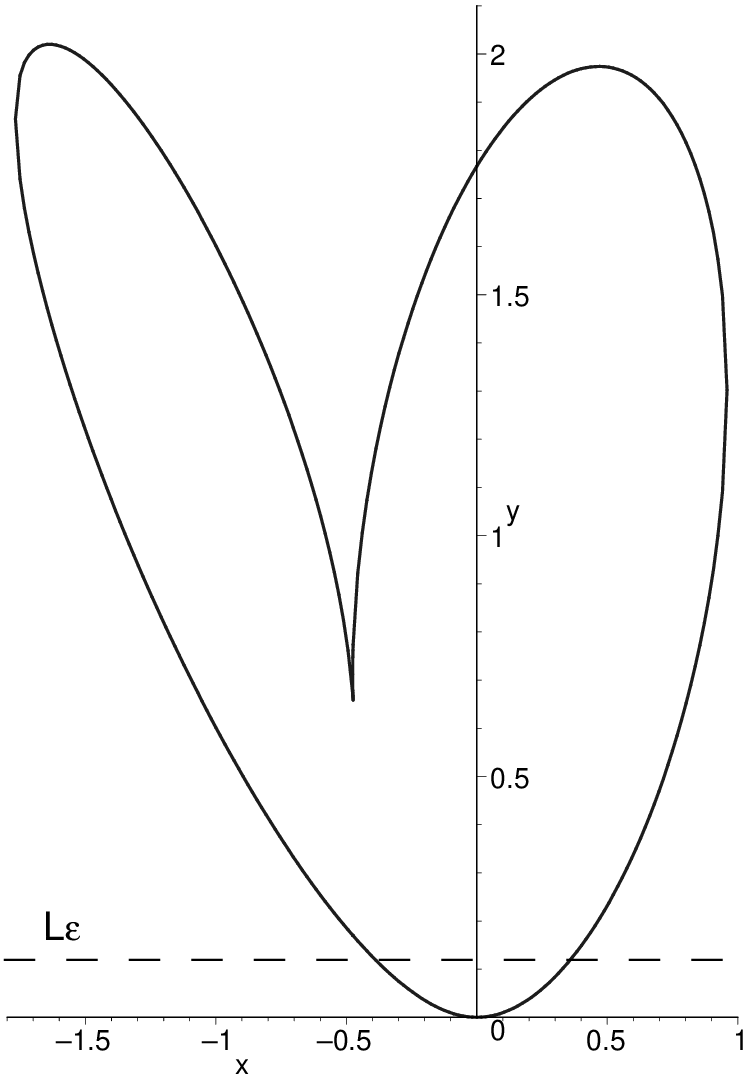}}
\put(0,120){\includegraphics{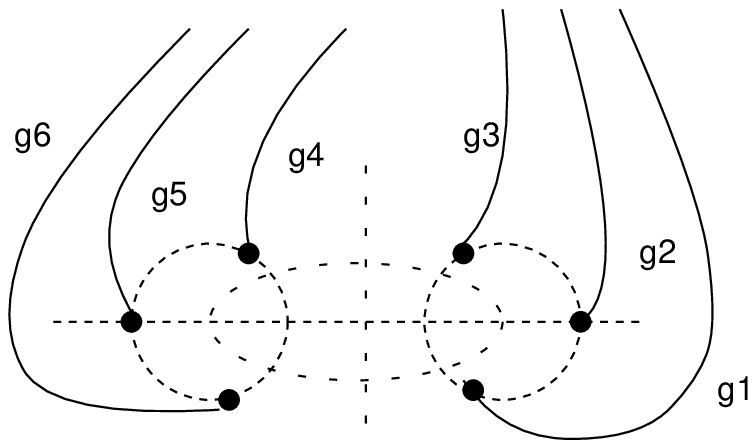}}
\end{picture}
\caption{Graph of $C_I$ and its generators at $y=\eps$}
\label{A2B(3,10)}
\end{figure}
\noindent
This curve $C_I$ has a $B_{3,10}$ singularity at $O$ and an
$A_2$ singularity at  $(-450/949,625/949)$.
We consider the pencils $L_t=\{y=t\}, t\in \bfC$.
% and take a generic pencil $L_\varepsilon, 0<\varepsilon\ll 1$. 
We take
  generators $g_1,g_2,\dots,g_6$ of the fundamental
group $\pi_1 (L_{\eta_0}-C_I)$ as in  Figure \ref{A2B(3,10)}
 where  $\eta_0 =\eps$.
Put $\omega_1:=g_3g_2g_1$ and $\omega_2:=g_6g_5g_4$. 
To  see the monodromy relation, we look at the Puiseux parametrization
at the origin. It is given by 
\[
 x=t^3+\frac{1}{2}t^7+\nnt, \quad y=t^6
\]
Therefore, at  $y=0$ we have  the following monodromy 
relations:
\begin{eqnarray} \label{A2B(3,10)-1}
g_1=g_6,\quad g_2=\omega_2 g_4 \omega_2^{-1},
\quad g_3=\omega_2 g_5 \omega_2^{-1}
\end{eqnarray} 
\begin{eqnarray} \label{A2B(3,10)-2}
g_4=(\omega_2\omega_1)g_3(\omega_2\omega_1)^{-1},\quad 
g_5=(\omega_2\omega_1^2)g_1(\omega_2\omega_1^2)^{-1},\quad 
g_6=(\omega_2\omega_1^2)g_2(\omega_2\omega_1^2)^{-1}
\end{eqnarray}
By taking the product of the relations in  (\ref{A2B(3,10)-1}),
% we have:
%\[
% g_3g_2g_1=\om_2g_5g_4\om_2^{-1}g_6
%\]
%Simplifying this,
 we get
\begin{eqnarray} \label{A2B(3,10)-good}
\om_1=\om_2
\end{eqnarray}
Then the big circle relation (=vanishing relation at infinity) 
$\om_2\om_1=e$ reduces to:
\begin{eqnarray} \label{A2B(3,10)-big}
\omega_1^2=e
\end{eqnarray} 
From  (\ref{A2B(3,10)-2}), (\ref{A2B(3,10)-good})
 and (\ref{A2B(3,10)-big})  we have 
\begin{eqnarray}\label{A2B(3,10)-2b}
g_4=g_3, \quad g_5=\omega_1 g_1 \omega_1^{-1}, \quad  
 g_6=\omega_1 g_2 \omega_1^{-1}
\end{eqnarray} 
From (\ref{A2B(3,10)-1}), (\ref{A2B(3,10)-good}) and (\ref{A2B(3,10)-2b})
we obtain $g_6=g_4=g_3=g_1$. Thus the generators are reduced to 
$g_1,g_2$. Rewriting   (\ref{A2B(3,10)-1}) as a relation for $g_1,g_2$, we get
%\begin{eqnarray*} 
%g_2=\om_2g_1\om_2^{-1}=\om_1g_1\om_1^{-1}=(g_1g_2g_1)g_1(g_1g_2g_1)\inv
%\end{eqnarray*}
%This implies
 the braid relation $\{g_1,g_2\}=e$.
On the other hand $e=\om_1^2=(g_1g_2g_1)^2$. 
Thus the fundamental group is isomorphic to $\Gr$ by Lemma \ref{max-equivalence}.

\vspace{.2cm}
\noindent
(II) {\bf Moduli space} $\cM(\{Sp_1,A_2\})$. 
Let us consider the following curve $C_{II}\in\cM(\{Sp_1,A_2\})  $
 which is defined by
\begin{eqnarray*}
&&C_{II}:f(x,y)=(y^2-y^3-x^3)^2-x^3y^3\\
&&\Delta_x(f)=-729y^{23}(y-1)^4(3y-4)^3, \quad P_y=\{0,1,4/3\}
\end{eqnarray*}

\begin{figure}[htb]
\setlength{\unitlength}{1bp}
\begin{picture}(255,130)(-100,0)
\put(-120,130){\includegraphics{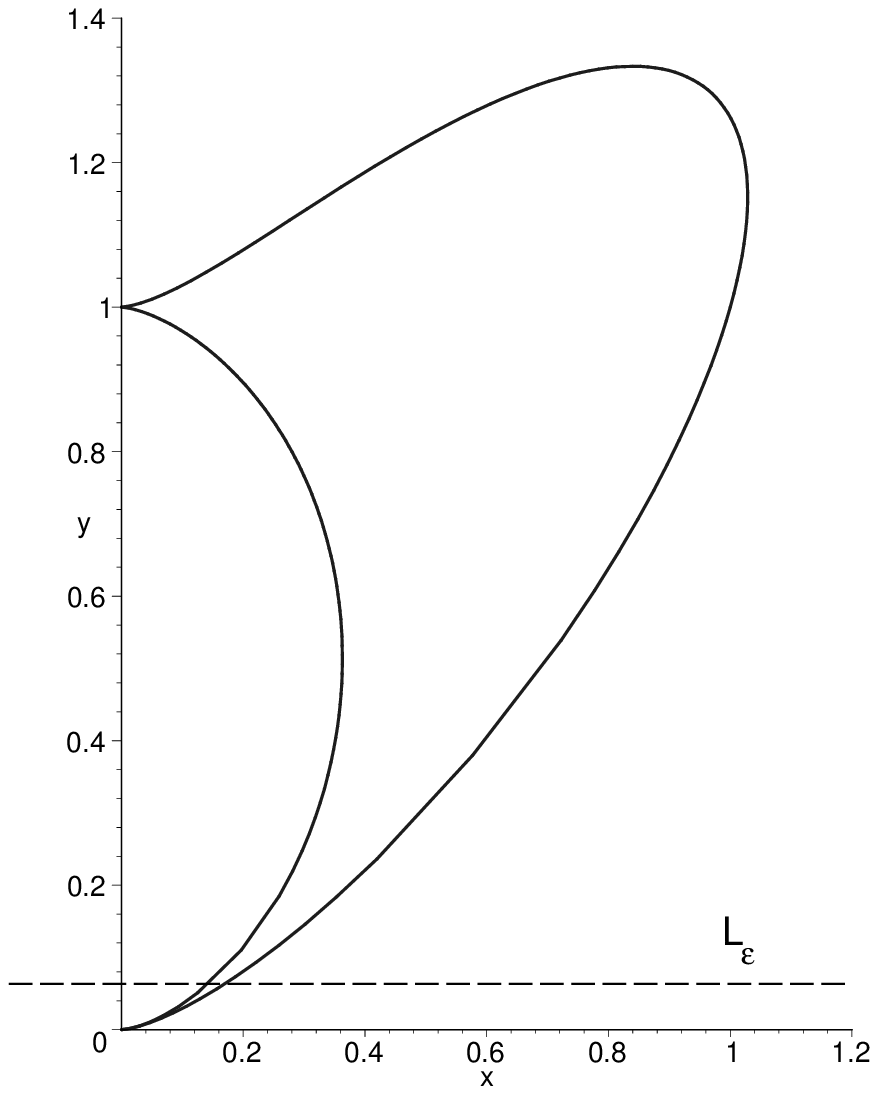}}
\put(40,130){\includegraphics{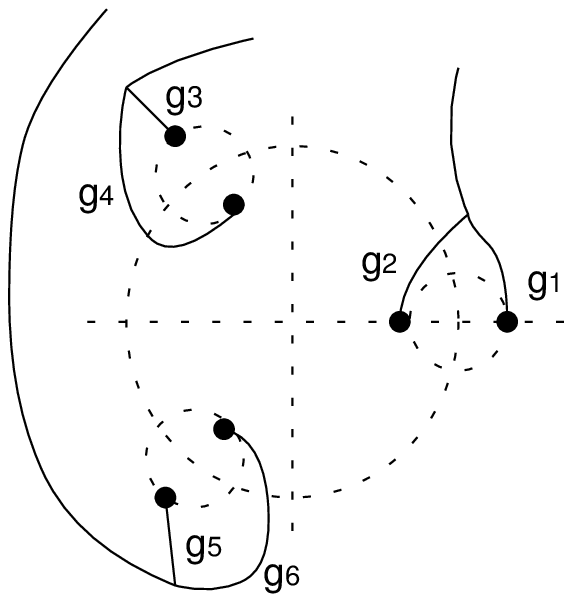}}
\end{picture}
\caption{Graph of $C_{II}$ and its generators at $y=\eps$}
\label{A2Sp1}
\end{figure}

This curve $C_{II}$ has  a $Sp_1$ singularity at the origin and 
an $A_2$ singularity at (0,1).
 We consider the pencils $L_t=\{y=t\}, t\in \bfC$.
% and  take a generic pencil $L_\varepsilon, 0<\varepsilon\ll 1$. 
We take
generators $g_1,g_2,\dots,g_6$ of the fundamental
group $\pi_1 (L_{\eta_0}-C_{II})$ as in  Figure 
\ref{A2Sp1}  where  $\eta_0 =\eps$. Put
$\om_1:=g_2g_1$, $\om_2:=g_4g_3$, $\om_3:=g_6g_5$ and
$\Om:=\om_3\om_2\om_1$. 

\noindent
To  see the monodromy relation at $y=0$ , we look at the Puiseux
parametrization at the origin. It is given by 
\[
 x=-t^4+\frac{1}{3}t^7+\nnt, \quad y=t^6
\]
The monodromy relations at $y=0$ are given by
\begin{eqnarray} \label{A2Sp1-1}
g_1=g_6,\quad g_2=\om_3 g_5\om_3\inv
\end{eqnarray}
\begin{eqnarray} \label{A2Sp1-2}
g_3=\Om g_2\Om\inv,\quad g_4=(\Om\om_1)g_1(\Om\om_1)\inv
\end{eqnarray}
\begin{eqnarray} \label{A2Sp1-3}
g_5=\Om g_4\Om\inv,\quad g_6=(\Om\om_2)g_3(\Om\om_2)\inv
\end{eqnarray}
Since the big circle relation is $\Om=e$, from the relations
(\ref{A2Sp1-1}), (\ref{A2Sp1-2}), (\ref{A2Sp1-3})  we get
\begin{eqnarray}
 g_6&=&g_1,\, g_3=g_2,\, g_5=g_4,\,\label{II-1}\\
 g_2&=&g_1g_4 g_1\inv,\,
g_4=g_2g_1 g_2\inv, \,
g_1=g_4g_2g_4\inv\label{II-2}
\end{eqnarray}
From (\ref{II-1}) and 
 the second relation of (\ref{II-2}),
we can reduce the generators to $g_1,g_2$ and 
 we obtain the braid relation
$\{g_1,g_2\}=e$ from (\ref{II-2}). On the other hand,
 the big
circle relation gives
   $(g_1g_2)^3=e$. Thus
 the fundamental group is isomorphic to $\Gr$
by Lemma \ref{max-equivalence}.

\vspace {0.2cm}
\noindent
(III) {\bf Moduli space} $\cM(\{B_{3,8},E_6\})$. 
Let us consider the curve which defined by 
\begin{eqnarray*}
&&C_{III}: f=((y-1)^2+x^2-1)^3+x^4y^2\\
&&\Delta_x(f)=-64y^{19}(y-2)^9(31y-54)^2, \quad P_y=\{0,54/31,2\}
\end{eqnarray*}
This curve $C_{III}$ has a $B_{3,8}$ singularity  at the origin and 
an $E_6$ singularity at (0,2).

\begin{figure}[htb]
\setlength{\unitlength}{1bp}
\begin{picture}(255,130)(-100,0)
\put(-150,130){\includegraphics{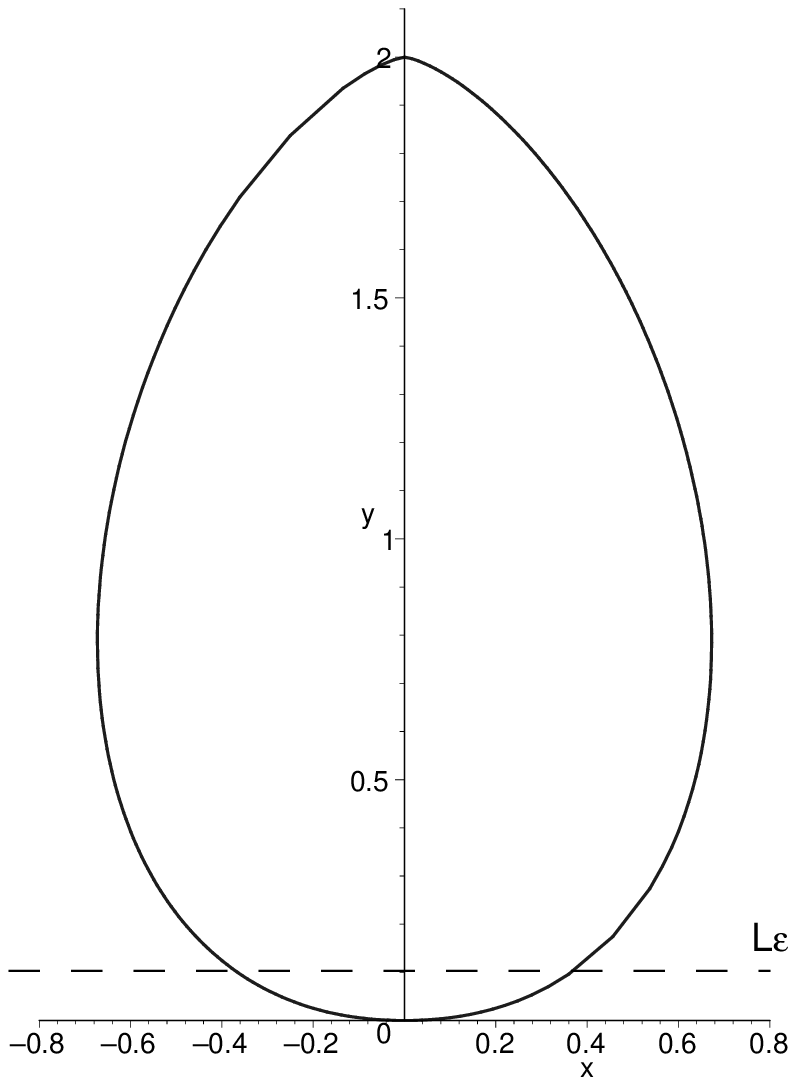}}
\put(0,130){\includegraphics{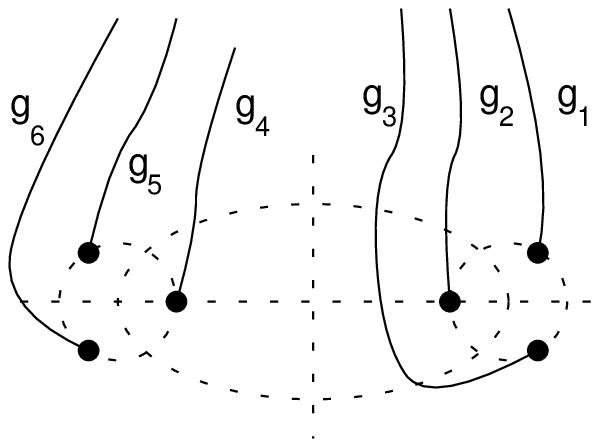}}
\end{picture}
\caption{Graph of $C_{III}$ and its generators at $y=\eps$}
\label{E6B(3,8)}
\end{figure}
\noindent
To  see the monodromy relation, we look at the Puiseux parametrization
at the origin. It is given by 
\[
 x=-\sqrt{2}t^3+\frac{\sqrt[6]{2}}{2}t^5+\nnt, \quad y=t^6
\]
%We consider the pencils $L_t=\{y=t\}, t\in \bfC$ and
%take a generic pencil $L_\varepsilon, 0<\varepsilon\ll 1$.
 We take 
generators $g_1,g_2,\dots,g_6$ of the fundamental
group $\pi_1 (L_{\eta_0}-C_{III})$ as in Figure
\ref{E6B(3,8)}
where  $\eta_0 =\eps$. 
Put $\omega_1:=g_3g_2g_1$ and $\omega_2:=g_6g_5g_4$. 
\newline
At the singular pencil line $y=0$, we have the  monodromy 
relations:
\begin{eqnarray} \label{E6B(3,8)-1}
g_1=\omega_2 g_4 \omega_2^{-1},\quad g_2=\omega_2 g_5 \omega_2^{-1},
\quad g_3=\omega_2 g_6 \omega_2^{-1}
\end{eqnarray} 
\begin{eqnarray} \label{E6B(3,8)-2}
g_4=\omega_2 g_3 \omega_2^{-1},\quad 
g_5=(\omega_2\omega_1)g_1(\omega_2\omega_1)^{-1},\quad 
g_6=(\omega_2\omega_1)g_2(\omega_2\omega_1)^{-1}
\end{eqnarray}
By taking the product of (\ref{E6B(3,8)-1}), we get $\om_2=\om_1$. Hence
the big circle relation is
\begin{eqnarray} \label{E6B(3,8)-big}
\omega_2\omega_1=\om_1^2=e
\end{eqnarray} 
From  (\ref{E6B(3,8)-2}) and (\ref{E6B(3,8)-big}) we get $g_5=g_1$ and 
$g_6=g_2$. Substituting these equalities to (\ref{E6B(3,8)-1}), we get
$g_1g_2g_1=g_2g_1g_4, ~ g_1g_4g_1=g_2g_1g_4$.
Thus we have $g_4=g_2$ and therefore the braid relation
$\{g_1,g_2\}=e$ follows.
Finally, from (\ref{E6B(3,8)-big})  we get $(g_2g_1)^3=e$.
 Thus
 the fundamental group is isomorphic to $\Gr$ by Lemma \ref{max-equivalence}.
 
\begin{Remark}  Recently,
Uluda\v{g} studied fundamental group of rational cuspidal curves.
Especially for sextics his results (see \cite[Theorem 3.1.5]{U}) is related to
our result. Namely, our cases $\{B_{3,10},A_2\}$, $\{B_{3,8},E_6\}$ and
$\{Sp_1,A_2\}$  correspond to his cases $\{[3_3],[2]\}$,
$\{[3_2,2],[3]\}$ and $\{[4,2_3],[2]\}$. He
showed the fundamental group is isomorphic to  $\bfZ_2*\bfZ_3$,
without assuming the curve to be a torus type.
\end{Remark}

\subsection{Reducible singularities} The last 4 moduli spaces contain 
reducible singularities, these singularities have 2 or 3 analytic
branches. The 
monodromy relation at these singularities is more complicate, but
the method we use is the same as before. For the computation of the
fundamental groups
for the curves in  $\cM(\{C_{9,9}\})$ and in  $\cM(\{C_{3,15}\})$, we use 
the existence of $\bfZ_2$ and $\bfZ_3$ actions
respectively on chosen curves.
We assume that $\eps$ is a positive, sufficiently small number as before. 

\vspace{.2 cm}
\noindent
(IV) {\bf Moduli space} $\cM(\{C_{3,7},A_8\})$. 
Let us consider the curve $C_{IV}\in \cM(\{C_{3,7},A_8\})$  defined by
\begin{eqnarray*}
C_{IV}:&f(x,y)&=(y(2y-2)-x^2)^3+\\
&&\left(\frac{46079}{54000}x^3+\frac{152279}{18000}x^2y+(\frac{311579}{18000}y^2-\frac{59}{10}y)x-\frac{2351327}{54000}y^3+\frac{178829}{4500}y^2\right)^2
\end{eqnarray*}
The discriminant polynomial is given by 
$\Delta_x(f)=cy^{16}(7y-6)^2(y-2)^9h_3(y)$
where $h_3(y)$ is a  polynomial of degree 3 which has 3 real
solutions, they are approximately 
$\beta_1\fallingdotseq 0.053$, $\beta_2\apeq 0.059$ and $\beta_3\apeq
0.831$.  This curve $C_{IV}$ has two singularities, 
a $C_{3,7}$ singularity at the origin and an $A_8$ singularity at
$(2,2)$.   Figure \ref{C(3,7)A8} shows the global
topological
(not numerical) situation
 of $C_{IV}$.

\begin{figure}[htb]
\setlength{\unitlength}{1bp}
\begin{picture}(255,130)(-100,0)
\put(-120,130){\includegraphics{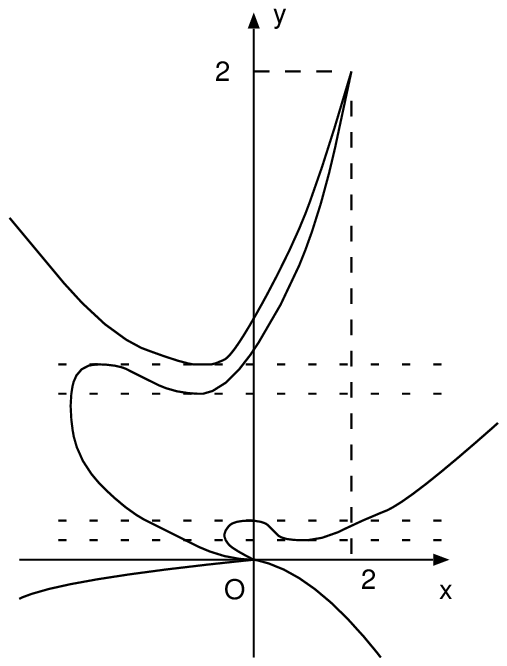}}
\put(30,130){\includegraphics{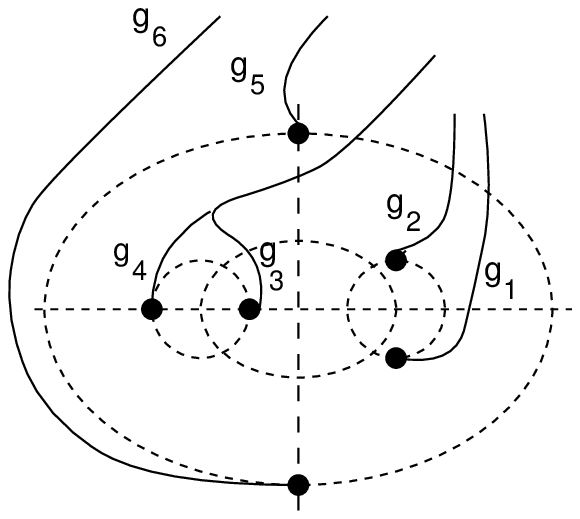}}
\end{picture}
\caption{Graph of $C_{IV}$ and its generators at $y=\eps$}
\label{C(3,7)A8}
\end{figure}
%We now going to show that the monodromy relations at $y=0$, $y=\beta_1$
%and $y=\beta_2$ are enough for compute the fundamental group. 
%
To see the monodromy relation at the origin, we look at the  Puiseux
expansion of $(C_{IV},O)\cong C_{3,7}$. It has two components, a smooth
component $L$ and a component $K$ of the (2,5)-cusp, which are 
 parametrized as follows.
\begin{eqnarray*}
&L:  &x=-\frac{60\sqrt{238}}{119}i\tau+\nnt, \quad y=\tau^2\\ 
&K:  &x=-\frac{30\sqrt{2}}{\sqrt{2581}}t^2+\frac{3481\sqrt[4]{1161450}}{38715}t^3+\nnt,\quad y=t^4
\end{eqnarray*}

We take  generators $g_1,g_2,\dots,g_6$ of the fundamental group 
$\pi_1(L_{\eta_0}-C_{IV})$ as in Figure \ref{C(3,7)A8}
where $\eta_0=\eps$ and $0<\varepsilon\ll 1$.
$g_5,g_6$ correspond to the points of $L$ and $g_1,\dots,g_4$
correspond to the points of $K$.
%
%We introduce two elements
Put  $\om_1:=g_2g_1$, $\om_2:=g_4g_3$.
The monodromy relations at the origin are given by:
\begin{eqnarray} \label{C3,7A8-O-L}
g_5=g_6 \quad \text{(relation for L)}
\end{eqnarray} 
\begin{eqnarray} \label{C3,7A8-O-K1}
g_1=g_5g_4g_5^{-1},\quad g_2=(g_5\om_2)g_3(g_5\om_2)^{-1} 
\end{eqnarray}
\begin{eqnarray} \label{C3,7A8-O-K2}
g_3=(g_5\om_2\om_1)g_1(g_5\om_2\om_1)^{-1},\quad 
g_4= (g_5\om_2\om_1)g_2(g_5\om_2\om_1)^{-1} 
\end{eqnarray}
Putting $h:=g_5=g_6$, the big circle relation is given by
\begin{eqnarray} \label{C3,7A8-big}
h^2\om_2\om_1=e
\end{eqnarray}
Then 
%we have $h\om_2\om_1=h^{-1}$ and 
(\ref{C3,7A8-O-K1}) and (\ref{C3,7A8-O-K2})  can be  rewritten as follows.
\begin{eqnarray} \label{C3,7A8-O-K1a}
g_1=hg_4h^{-1},\quad g_2=(h\om_2)g_3(h\om_2)^{-1} 
\end{eqnarray}
\begin{eqnarray} \label{C3,7A8-O-K2a}
g_3=h^{-1}g_1h,\quad 
g_4= h^{-1}g_2h 
\end{eqnarray}
From  (\ref{C3,7A8-O-K1a}) and (\ref{C3,7A8-O-K2a}) we get 
$\rho:=g_2=g_1$, $\xi:=g_4=g_3$. 
%So we can reduce to \rho,\xi$.
 We rewrite (\ref{C3,7A8-O-K1a}) and 
(\ref{C3,7A8-big}) as follows.
\begin{eqnarray} \label{C3,7A8-O-K1b}
\rho=h\xi h^{-1} 
\end{eqnarray}
\begin{eqnarray} \label{C3,7A8-O-K2b}
h^2\xi^2\rho^2=e
\end{eqnarray}
To read the monodromy at $y=\beta_1$ we need to know how the generators
 move, when the pencil line $y=\eta$ move from
 $\eta=\eps\to\beta_1^{-}:= \beta_1-\eps$.
We can show  see
the generators are deformed
as in Figure \ref{C(3,7)A8-mov}, using a similar argument as in 
\cite{Oka-sextics}.

\begin{figure}[htb]
\setlength{\unitlength}{1bp}
\begin{picture}(255,100)(-100,0)
\put(-120,100){\includegraphics{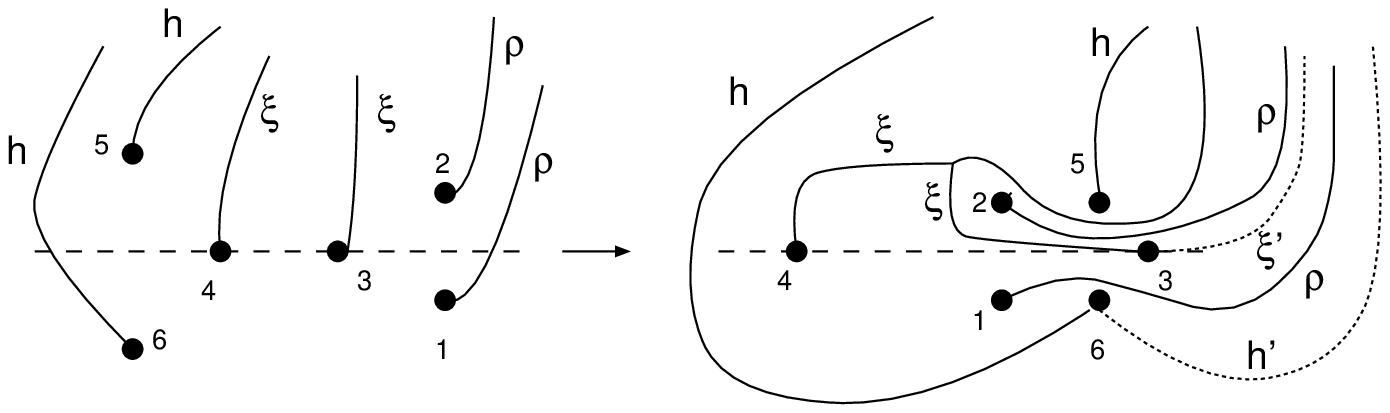}}
\end{picture}
\caption{The movement of the generators $y:\eps\to\beta_1^{-}$}
\label{C(3,7)A8-mov}
\end{figure}
Now  the monodromy relations at $y=\beta_1$ and $y=\beta_2$
 are simple tangent relations, which are given by
\begin{eqnarray} \label{C3,7A8-flex}
h=h'=\xi' 
\end{eqnarray}
where $h'=h$ (because of the big circle relation) and 
$\xi'=\rho\inv\xi\rho$.
Thus $h=\rho\inv\xi\rho$.
Thus we can take $\rho,\xi$ as the generators.
 Now apply to (\ref{C3,7A8-O-K1b}) and
(\ref{C3,7A8-O-K2b}), we obtain
$\rho\xi\rho=\xi\rho\xi$ and  
$(\rho\xi\rho)^2=e$. Thus the fundamental group is isomorphic to
$\Gr$
by Lemma \ref{max-equivalence}.

\vspace{.2cm}
\noindent
(V) {\bf Moduli space} $\cM(\{C_{9,9}\})$. We will use
% the symmetric curve, and apply
 the technique in \cite{Two}. Let us consider the affine curve
\[C_g: g=(y-x^2)^3+\left(2y-2x^2+\frac{32}{27}x^3\right)^2=0\]
and let $f(x,y):=g(x,y^2)$ and 
 let $C_{V}$ be the curve defined by $ f(x,y)=0$.
We can easily observe  that 
$(C_g,O)\cong A_8$ and $(C_V,O)\cong C_{9,9}$.
 The discriminant polynomial
of $g$ is given by
\[\Delta_x g=c y^{11}(y+4)(295y^3+208y^2-3456y+729)\]
Thus $\Delta_x g$ has 5 real solutions, they are approximately  $\gamma_0=0$,
$\gamma_1\apeq 0.21$, $\gamma_2\apeq 2.96$, $\gamma_3\apeq -3.88$ and
$\gamma_4=-4$. 
We denote
$\delta_i=\sqrt{\gamma_i}$ ($i=1,2,3,4$) and $\delta_0:=\gamma_0$. 
An easy observation that 
$P_y(f)=\{\delta_0,\pm\delta_i\}$, for $i=1,2,3,4$. 
%Positions of $\gamma_i,\pm\delta_i$ in a complex plane
%are showed in
See  Figure  \ref{C9,9-gen-positions}.
\begin{figure}[htb]
\setlength{\unitlength}{1bp}
\begin{picture}(255,130)(-100,0)
\put(-100,130){\includegraphics{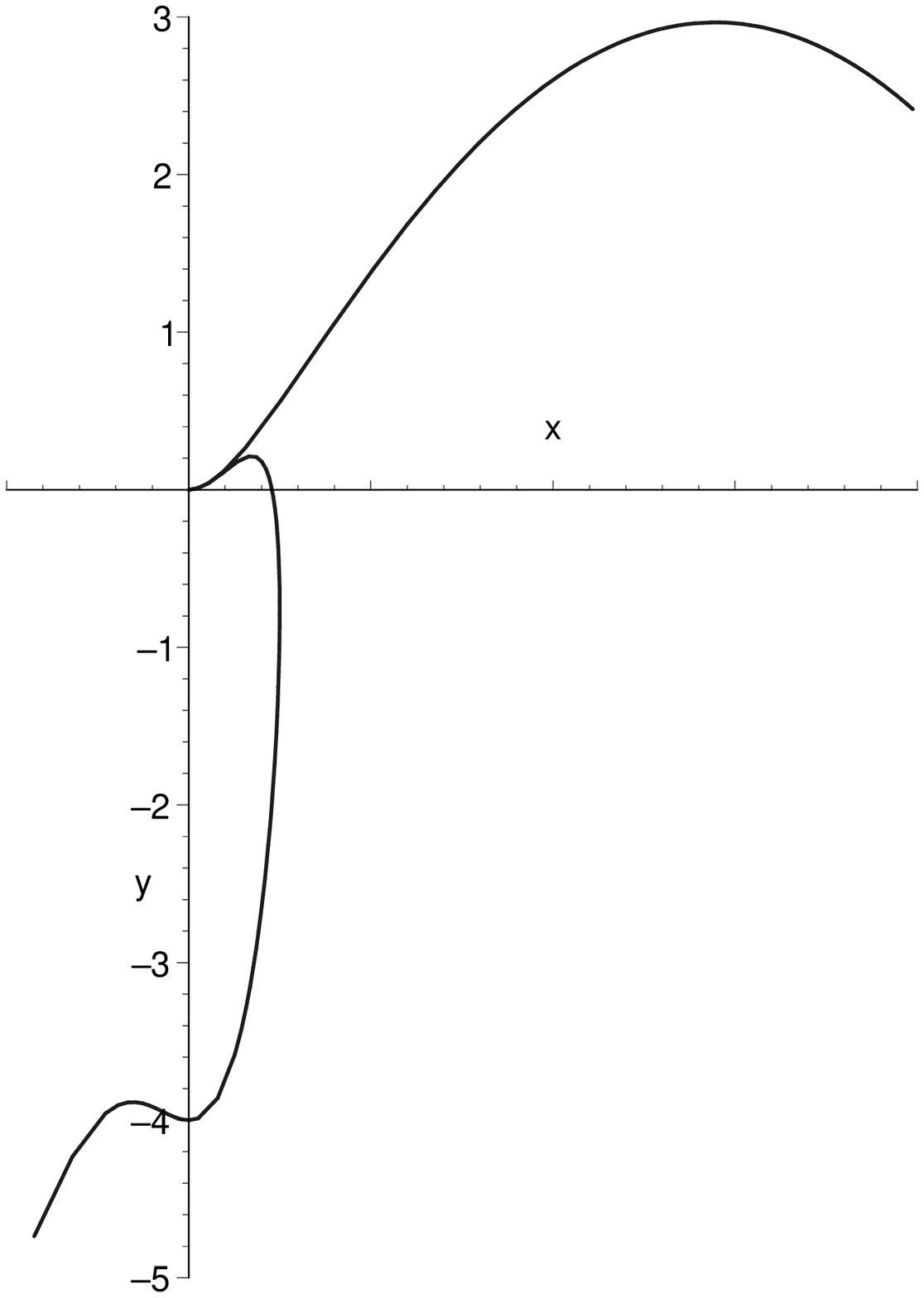}}
\put(40,130) {\includegraphics{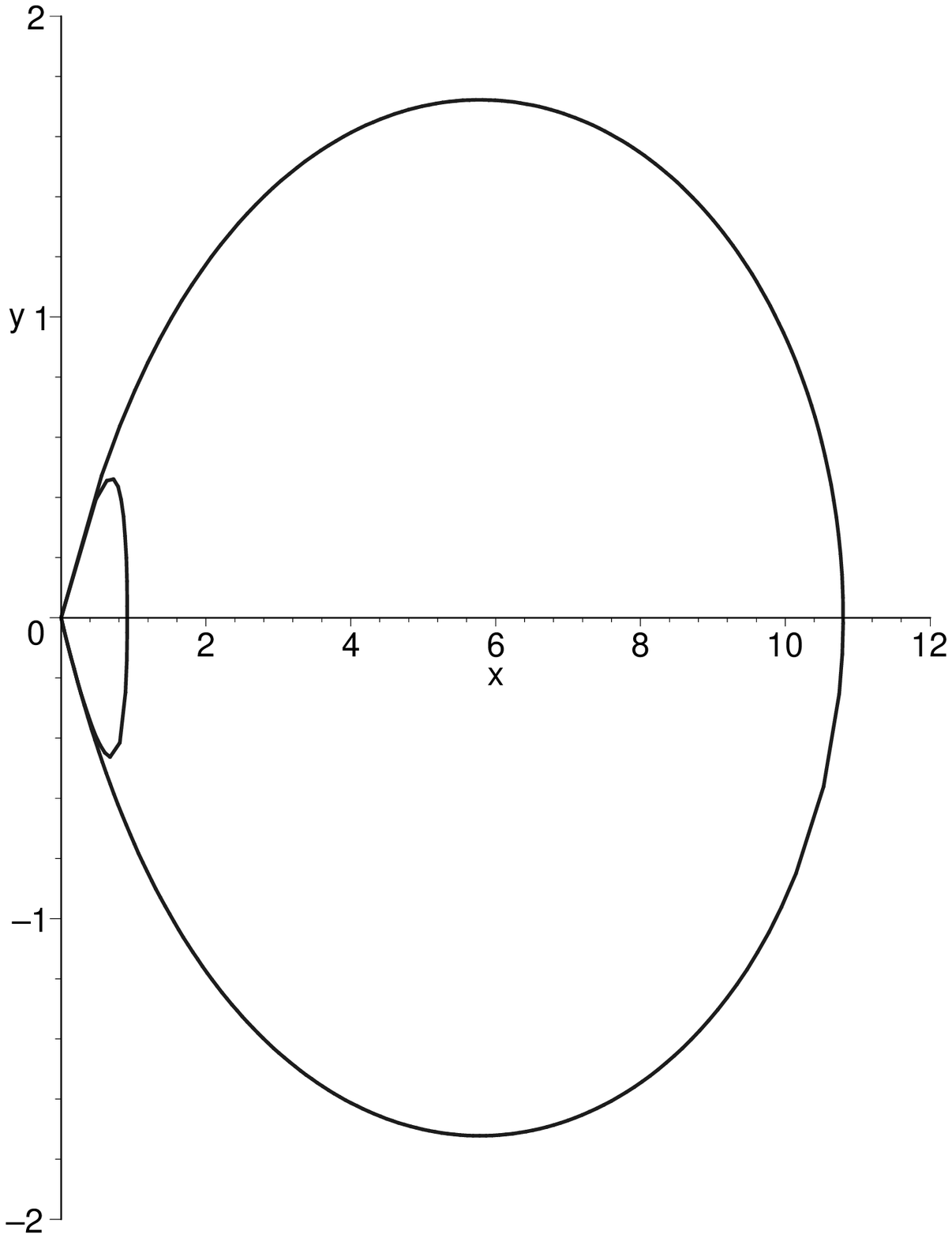}}
\end{picture}
\caption{Graph of $C_g$ and $C_{V}$}
\label{C9,9}
\end{figure}
Since $C_{V}$ is a double covering of $C_g$ along the $x$-axis, we are
able to read the monodromy relations of $C_{V}$ via $C_g$. 
We consider the pencils $L_t=\{y=t\}, t\in \bfC$.
% Let take a generic 
%pencil $L_{\varepsilon}$, $0<\varepsilon\ll 1$. 

\begin{figure}[htb]
\setlength{\unitlength}{1bp}
\begin{picture}(255,130)(-100,0)
\put(-100,130){\includegraphics{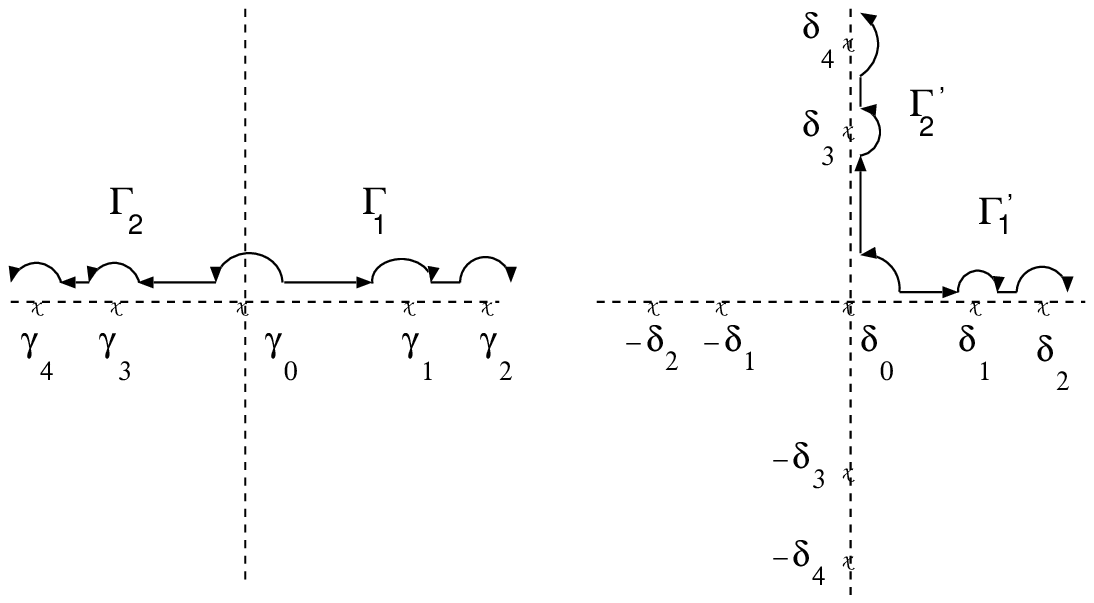}}
\end{picture}
\caption{Position of singular pencil lines of $C_g$ and $C_{V}$}
\label{C9,9-gen-positions}
\end{figure}

The monodromy relations at $y=\delta_i$
($i=0,1,\dots,4$) are enough to compute the fundamental group.
We first observe that
the monodromy relation at $y=\delta_i$ for the curve $C_{V}$ is nothing 
but   the monodromy relation  for the curve $C_g$ at $y=\gamma_i$ ($i\ne
0$). 
\begin{figure}[htb]
\setlength{\unitlength}{1bp}
\begin{picture}(255,130)(-100,0)
\put(-100,130){\includegraphics{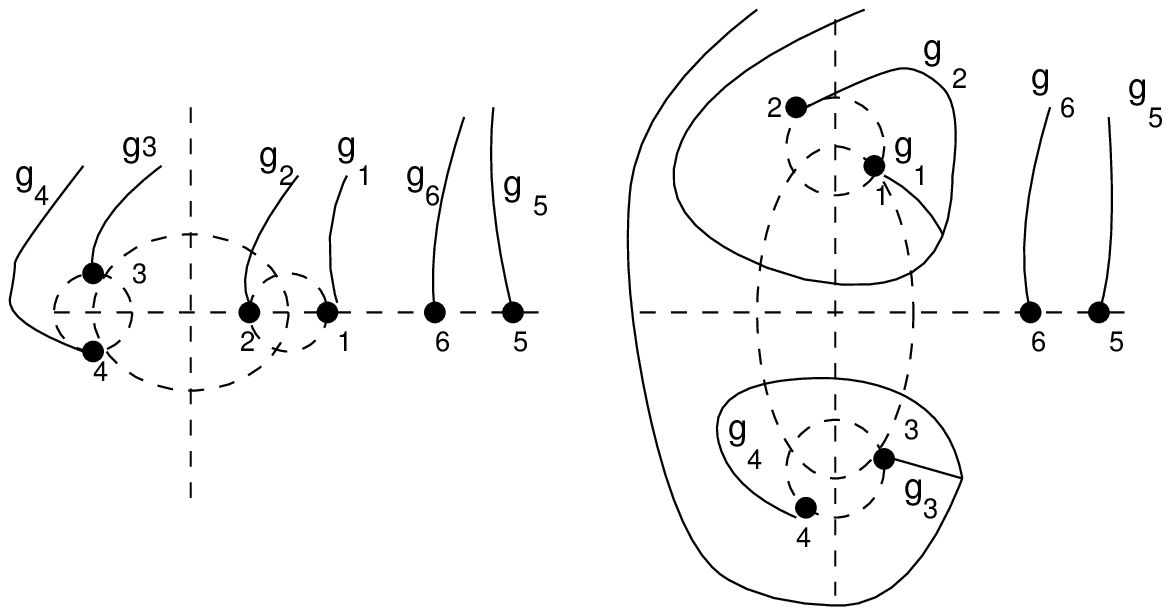}}
\end{picture}
\caption{Generators at $y=\eps$ and $y=-\eps$}
\label{C9,9-gen}
\end{figure}

To  see the monodromy relation at $y=0$, we look at the Puiseux parametrization
of $C_g$ at the origin. It is given by
\[
 x=\phi(t),\quad \phi(t):=t^2+\frac{16\sqrt{3}}{243}t^7+\nnt, \quad y=t^4
\]
We take  generators $g_1,g_2,\dots,g_6$ of the fundamental group 
$\pi_1(L_{\eps}-C_{V})$ as in Figure \ref{C9,9-gen}.
Put $\om_1:=g_2g_1$ and  $\om_2:=g_4g_3$.
%and $\Om=\om_2\om_1$.
%
Note that $C_{V} $ is parametrized  at the origin  by 
$x(t)=\phi(t),~y(t)=t^2$ and $(C_V,O)$ has two irreducible components. Thus
the monodromy relations at $y=0$ for $C_V$ are given by
\begin{eqnarray} \label{C9,9-O1}
g_1=\om_2\om_1^3g_2(\om_2\om_1^3)^{-1},\quad  g_2=\om_2\om_1^4g_1(\om_2\om_1^4)^{-1}
\end{eqnarray}
\begin{eqnarray} \label{C9,9-O2}
g_3=\om_2\om_1\om_2^2g_4(\om_2\om_1\om_2^2)^{-1},
\quad  g_4=\om_2\om_1\om_2^3 g_3(\om_2\om_1\om_2^3)^{-1}
\end{eqnarray}
The tangent relation at $y=\delta_1$:
\begin{eqnarray} \label{C9,9-P1}
g_6=g_1
\end{eqnarray}
The tangent relation at $y=\delta_2$:
\begin{eqnarray} \label{C9,9-P2}
g_5=g_1^{-1}g_2g_1
\end{eqnarray}
From (\ref{C9,9-P1}) and (\ref{C9,9-P2}), the big circle
relation  $\om_2\om_1g_6g_5=e$ and (\ref{C9,9-O1}) reduce to:
\begin{eqnarray} \label{C9,9-big}
\om_2\om_1^2=e,\quad
g_1g_2g_1=g_2g_1g_2
\end{eqnarray}
Also  the relation (\ref{C9,9-O2})  can be written as 
\begin{eqnarray} \label{C9,9-O4}
\om_1^5g_3=g_4\om_1^5
\end{eqnarray}
The tangent relation at $y=\delta_3$ (when 2 points which have indices 2 and 4
in Figure \ref{C9,9-gen} coincide):
\begin{eqnarray} \label{C9,9-P3}
g_3^{-1}g_4g_3=g_1^{-1}g_2g_1
\end{eqnarray}
\noindent
From (\ref{C9,9-P3}) and (\ref{C9,9-O4}) we can have
\begin{eqnarray} \label{C9,9-g3g4}
g_3=\om_1^{-3}g_1,\quad g_4=\om_1^2g_1\om_1^{-5}
\end{eqnarray}
Finally we rewrite (\ref{C9,9-big}) and obtain
\begin{eqnarray} \label{C9,9-last}
\om_1^2g_1\om_1^{-8}g_1\om_1^2=e
\end{eqnarray}
Using  the second relation of (\ref{C9,9-big}) we get $(g_2g_1)^3=e$ and we get a canonical
surjection from $\bfZ_2*\bfZ_3$ to
$\pi_1(\bfP^2-C_V)$.
Thus the fundamental group is isomorphic to $\Gr$ by Lemma \ref{max-equivalence}.

\vspace{.2cm}
\noindent
(VI)  {\bf Moduli space   $\cM(\{C_{3,15}\})$}.
We consider the
following family of tame sextics curve  $C_{VI,t}\in\cM(\{C_{3,15}\}) $ with a
$C_{3,15}$-singularity
 at the origin
 which is define by
\begin{eqnarray*}
C_{VI,t}:
f(x,y)&=&f_2(x,y)^3+f_3(x,y,t)^2=0\\
f_2(x,y)&=&y-x^2 ,\quad
f_3(x,y,t) = t yx+2 y^3-t x^3,\quad 1\le t\le \frac{3\sqrt{2}}{4}
\end{eqnarray*}
This family enjoys the following  properties:

(1) For  any $t$, $C_{VI,t}$ has a unique $C_{3,15}$ singularity at the origin $O$.
For $t\ne1$,
$C_{VI,t}$ is irreducible but
$C_{VI,1}$ is reducible and it consists of a  line $y=0$ and a quintic with
an  $A_{12}$-singularity at
the origin. However the local singularity at the origin is still 
$C_{3,15}$.

(2) $C_{VI,t}$ is stable under the $\bfZ_3$ action which is induced by
$(x,y)\mapsto (x\alpha, y\alpha^2)$ where $\alpha:=\exp(2\pi i/3)$.

For a practical computation, we take
$t=t_0$ with $t_0:=3\sqrt{2}/4$. 
\begin{figure}[htb]
\setlength{\unitlength}{1bp}
% \begin{picture}(255,255)(0,-6)
 \begin{picture}(255,130)(-100,0)
\put(-130,130){\includegraphics{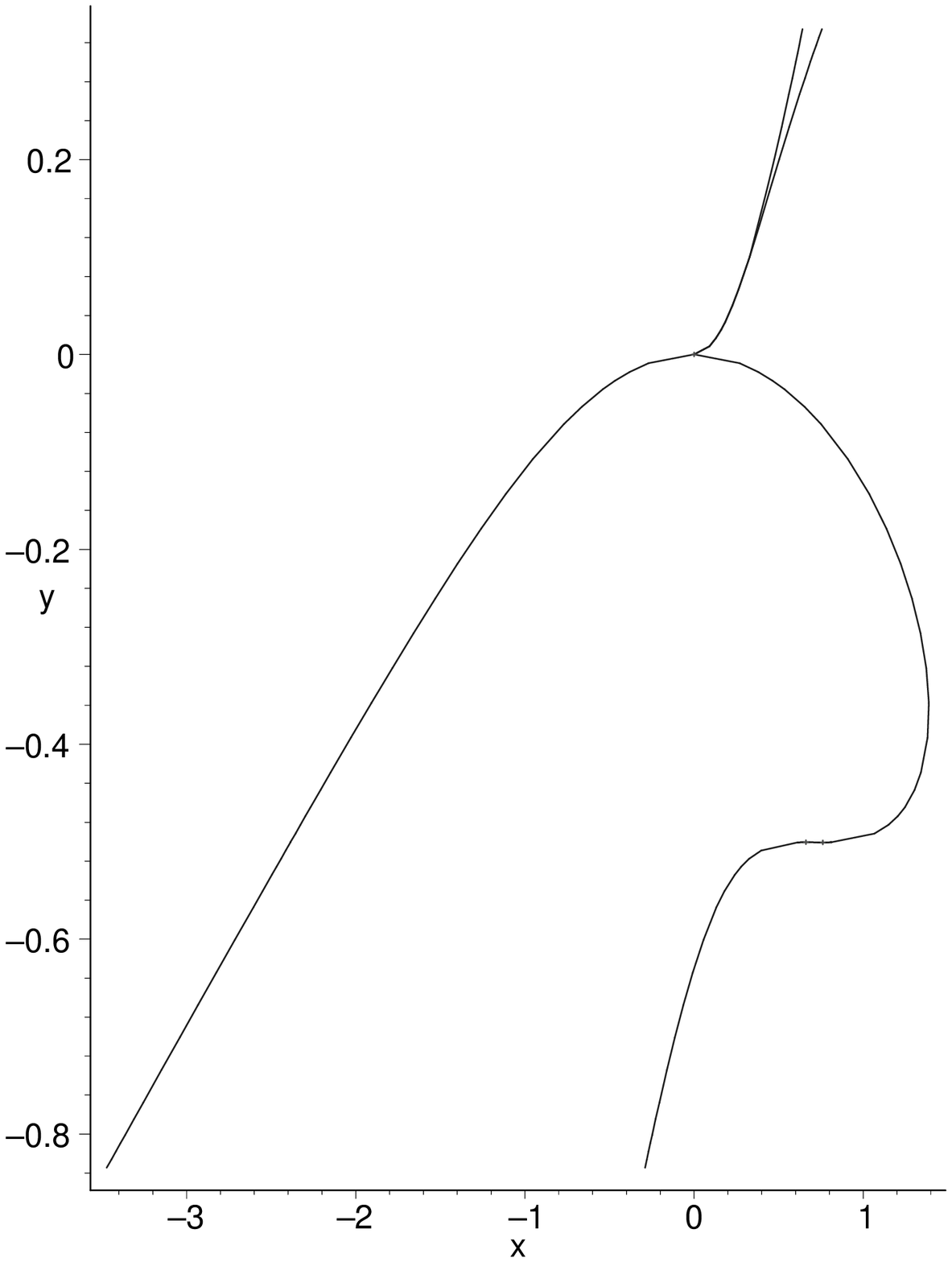}}
\put(10,130){\includegraphics{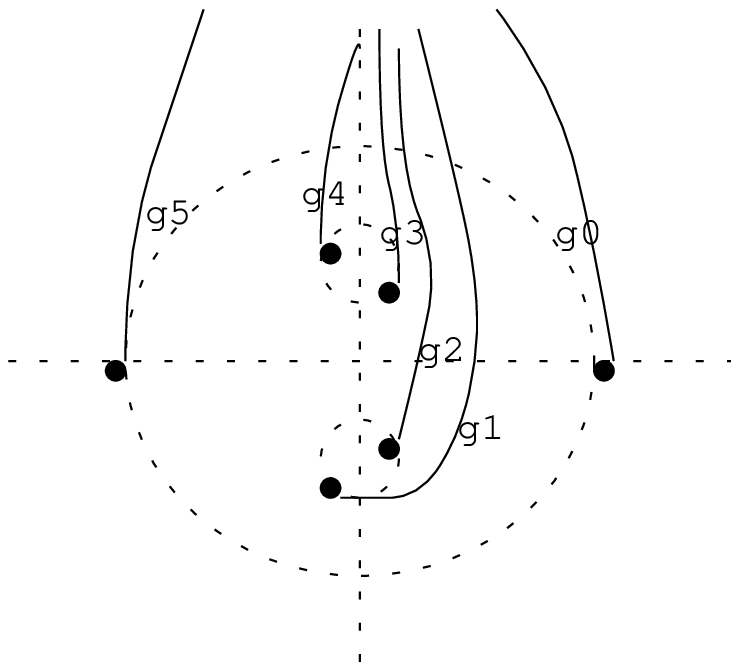}}
\end{picture}
%\vspace{2cm}
\caption{Graph of $C_{VI,t_0}$ and its generators at $y=-\eps$}
\label{C3,15graph}
\end{figure}
\noindent
We consider the pencil $L_t= \{y=t\},~t\in \bfC$.
The discriminant polynomial
$\De_x(f)$ is given as
\begin{eqnarray*}
\De_x(f)=&11664 y^{24}(1+(2y)^3)^2
%11664 y^{24}(1+2y)^2 (1-2y+4y^2)^2
\end{eqnarray*}
Thus the singular singular pencil lines $ y=-\alpha^j/2,j=0,1,2$ which are the
tangent lines at the flexes.
% and  two of them are  non-real singular pencil lines.

To see the monodromy relations at $y=0$, we use at the Puiseux
parametrizations of $C_{VI,t_0}$ at the origin.
It has 
a smooth component $L$ and a component $K$ of the $(2,13)$-cusp.
They have the following parametrizations by Lemma \ref{inverse-Puiseux}:
\begin{eqnarray*}
&L:  &x=-2\sqrt{2} i \tau+\text{(higher terms)}, \quad y=\tau^2\\ 
&K: 
&x=s^2+\frac{2\sqrt{2}}{3}s^8+\frac{11\sqrt{3\sqrt{2}}}{27} s^{11}
+\text{(higher terms)},\quad y=s^4
\end{eqnarray*}
We take  generators $g_0,g_1,\dots,g_5$ of the fundamental
group 
$\pi_1 (L_{\eta_0}-C)$ with $\eta_)=-\eps$ as in Figure \ref{C3,15graph}.
% where  $\eta_0 =-\eps$ and $\eps>0$  are assumed to be
%sufficiently small as before.
 Note that 
$g_0,g_5$ correspond to the points of $L$ and $g_1,\dots, g_4$ correspond  to
the points of $K$. 
We can show $\pi_1(\bfP^2-C_{VI,t_0})\cong \bfZ_2*\bfZ_3$ by a similar computation as
before.

Instead of giving a boring computation, we give a simpler proof
 by considering the degenerated curve $C_{VI}:=C_{VI,1}$.
As  $C_{VI}$ is reducible and  has the line $y=0$ as a component,
we take $y=0$ as the line at infinity.
Then \[
\pi_1(\bfP^2-C_{VI})\cong \pi_1(\bfC_y^2-Q)\]
where $\bfC_y^2$ is the affine chart with $(x,z)$ as coordinates and 
$Q$ is a quintic  which is defined by 
\begin{eqnarray}
Q:\quad z^3-2 z^2 x^2+z x^4+4 x z+4-4 x^3=0
\end{eqnarray}
This quintic is a rational curve with an $A_{12}$ singularity at the infinity.
We consider the pencil $x=\eta, \eta\in \bfC$. It has three simple tangents defined by
$32 x^3-432=0$.
We take $O=(0,0)$ to be the base point of the fundamental group and 
take three generators $g_0,g_1,g_3$ on the pencil line
$x=0$ as in Figure \ref{Quintics}. 
As $Q$ has also $\bfZ_3$ symmetry defined by 
$(x,z)\mapsto (x\alpha,z\alpha^2)$,
the monodromy relations are given by 
\[ g_0=g_1g_2g_1\inv,\quad  g_1=g_2g_0g_2\inv,\quad g_2=g_0g_1g_0\inv\]
We can immediately see that 
\[\pi_1(\bfC^2-Q)\cong \langle g_0,g_1;g_0g_1g_0=g_1g_0g_1\rangle
\cong B(3)\]
As  $D(B(3))\cong F(2)$, we get the surjective
homomorphism:
\[
F(2)\cong D(\pi_1(\bfC^2-Q))= D(\pi_1(\bfP^2-C_{VI}))
\to D(\pi_1(\bfP^2-C_{VI,t_0}))
\]
this implies that 
$\pi_1(\bfP^2-C_{VI,t_0}))\cong \bfZ_2*\bfZ_3$ by Lemma \ref{max-equivalence}.
\begin{figure}[htb]
\setlength{\unitlength}{1bp}
 \begin{picture}(255,120)(-100,0)
\put(-120,120){\includegraphics{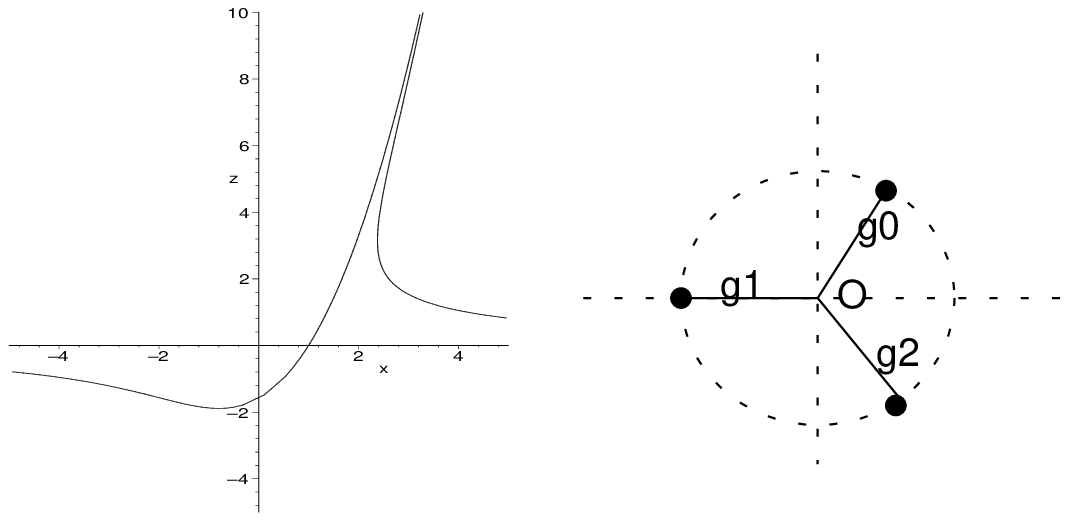}}
\end{picture}
\vspace{.4cm}
\caption{Quintics and generators at $x=0$}
\label{Quintics}
\end{figure}

\vspace {0.2cm}
\noindent
(VII) {\bf Moduli space} $\cM(\{C_{3,8},3A_2\})$. 
Let us consider the curve $C_{VII}\in\cM(\{C_{3,8},3A_2\}) $
 which is defined by
\begin{eqnarray*}
C_{VII}: f(x,y)=
 (y-x^2)^3+
\left( \frac{\sqrt[4]{27}}{423} y^3+(\frac{\sqrt[4]{3}}{4}+
\frac{1}{72} x)y^2+(\frac{3\sqrt{3}}{2} x+ \frac{\sqrt[4]{3}}{2} x^2)y
\right)^2
\end{eqnarray*}
\[
\Delta_x(f)=
c y^{17} (239 y^2+10800\sqrt{3} y -139968)(5y+108\sqrt{3})^2(y-4\sqrt{3})^3
(y+108\sqrt{3})
^6
\]
Note that 
$P_y=\{\al_1,\al_2,\al_3,\al_4,\al_5,\al_6\}$
where 
$\al_1:=-108\sqrt{3}$,
$\alpha_2:=-\frac{5400\sqrt{3}}{239}-\frac{7776\sqrt{2}}{239}$,
$\al_3:=-\frac{108\sqrt{3}}{5}$,
$\al_4:=0$, $\al_5:=-\frac{5400\sqrt{3}}{239}+\frac{7776\sqrt{2}}{239}$ and 
$\alpha_6:=4\sqrt{3}$. We observe that 
$\al_1<\al_2<\al_3<\al_4=0<\al_5<\al_6$.
This curve $C_{VII}$ has a $C_{3,8}$ singularity  at the origin and three
$A_2$ singularities, namely  one $A_2$ in the level $y=\al_6$ and
two  $A_2$ singularities in the level $y=\al_1$. The other singular
pencils are $y=\al_2,\al_5$ (simple tangents) and
$y=\al_3$ (double tangent).
\begin{figure}[htb]
\setlength{\unitlength}{1bp}
\begin{picture}(255,130)(-100,0)
\put(-120,130){\includegraphics{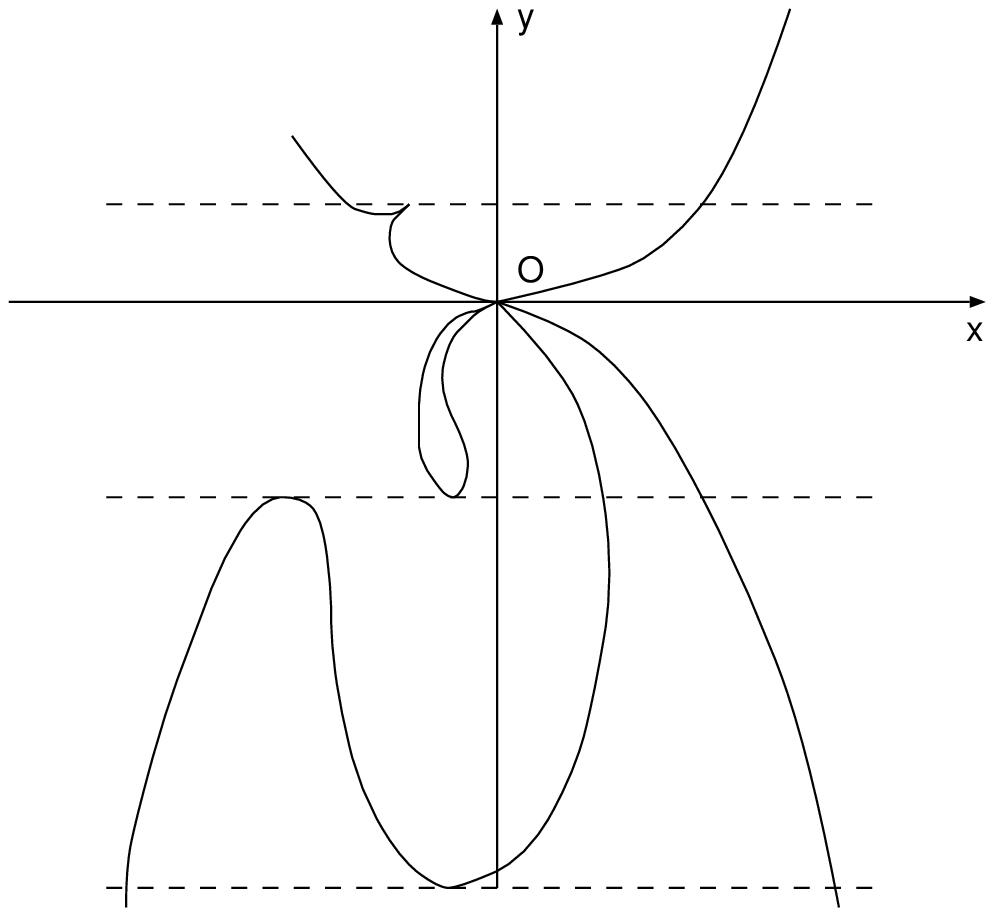}}
\put(20,130){\includegraphics{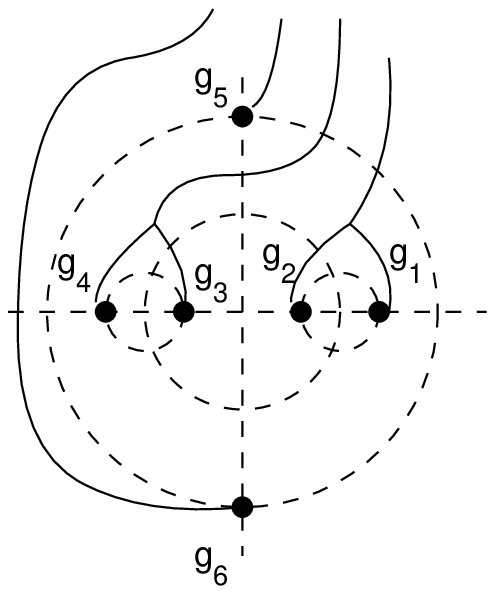}}
\end{picture}
\caption{Graph of $C_{VII}$ and its generators at $y=-\eps$}
\label{C(3,8)3A2}
\end{figure}
We are going to show that the monodromy relations at $y=0$,
$y=-\frac{108\sqrt{3}}{5}$ and $y=\alpha_2$ are enough to
compute the fundamental group. 
\newline
We take  generators $g_1,g_2,\dots,g_6$ of the fundamental group 
$\pi_1(L_{\eta_0}-C_{VII})$  as in Figure \ref{C(3,8)3A2}
where
$\eta_0=-\eps$.
% and $0<\varepsilon\ll 1$.
Since  $y=\al_3$ is
a double tangent, 
we obtain the relations
%that 6 generators contain: 
$g_5=g_6$ and $g_3=g_4$. Put
$\rho:=g_5=g_6,\quad \xi:=g_3=g_4$.
%\newline
The tangent relation at $y=\alpha_2$ is given by  
\begin{eqnarray} \label{3A2C3,8-simpletangent}
g_2=\xi\inv \rho\xi
\end{eqnarray}
To see the monodromy relation at the origin, we look at the  Puiseux
expansion of $(C_{VII},O)\cong C_{3,8}$. It has three components, a smooth
component $L$ and two smooth
components $K_1, K_2$.
% come from the (2,6)-cusp, 
Moreover we can compute explicitly 
their parametrization, 
\begin{eqnarray*}
&L:  &x=2t^2+\nnt, \quad y=\tau^2\\ 
&K_{i}:  &x=-\frac{i\sqrt{2}}{2}t^2+ 
z_{i}t^4+\nnt,\quad y=t^4,\quad z_{i}\ne 0,~i=1,2
\end{eqnarray*}
Position of the generators at $y=-\eps$ are showed in Figure \ref{C(3,8)3A2}.
Thus the monodromy relations at $y=0$ are given by:
\begin{eqnarray} \label{3A2C3,8-O}
g_2=g_1=\rho\xi\rho\inv
\end{eqnarray}
Finally, using the big loop relation $\rho^2\xi^2g_2g_1=e$, and together with  
(\ref{3A2C3,8-simpletangent}) and (\ref{3A2C3,8-O}), we obtain
 $\{\rho,\xi\}=e$ and $(\xi\rho\xi)^2=e$.
Thus the fundamental group is isomorphic to $\Gr$ by 
Lemma \ref{max-equivalence}.

%\vspace {0.2cm}
%\noindent
\subsection {Exceptional moduli space $\cM(\{C_{3,9},3A_2\})$}
%To compute the fundamental group
%$\pi_1(\bfP^2-C)$ for a curve $C\in \cM(\{C_{3,9},3A_2\})$,
We take the curve $C_{\text{ex}}\in\cM(\{C_{3,9},3A_2\}) $ which is defined by
\begin{eqnarray*}
C_{\text{ex}}: f(x,y)=(y-x^2)^3+
\left( y^2+\frac{4}{3} y^3+\frac{3\sqrt{3}}{2}xy+\frac{2\sqrt{3}}{3}xy^2+2x^2y
\right)^2
\end{eqnarray*}
\[
\Delta_x(f)=\frac{8192}{19683}y^{18}(2y+3)^2(4y-3)^4(4y+9)^6
\]
Thus the singular pencil lines correspond to $P_y=\{-9/4,-3/2,0,3/4\}$.
This curve $C_{\text{ex}}$ has a $C_{3,9}$ singularity  at the origin and three
$A_2$ singularities (where one $A_2$ in the level $y=3/4$,
two  $A_2$ in the level $y=-9/4$).

\begin{figure}[htb]
\setlength{\unitlength}{1bp}
\begin{picture}(255,130)(-100,0)
\put(-120,130){\includegraphics{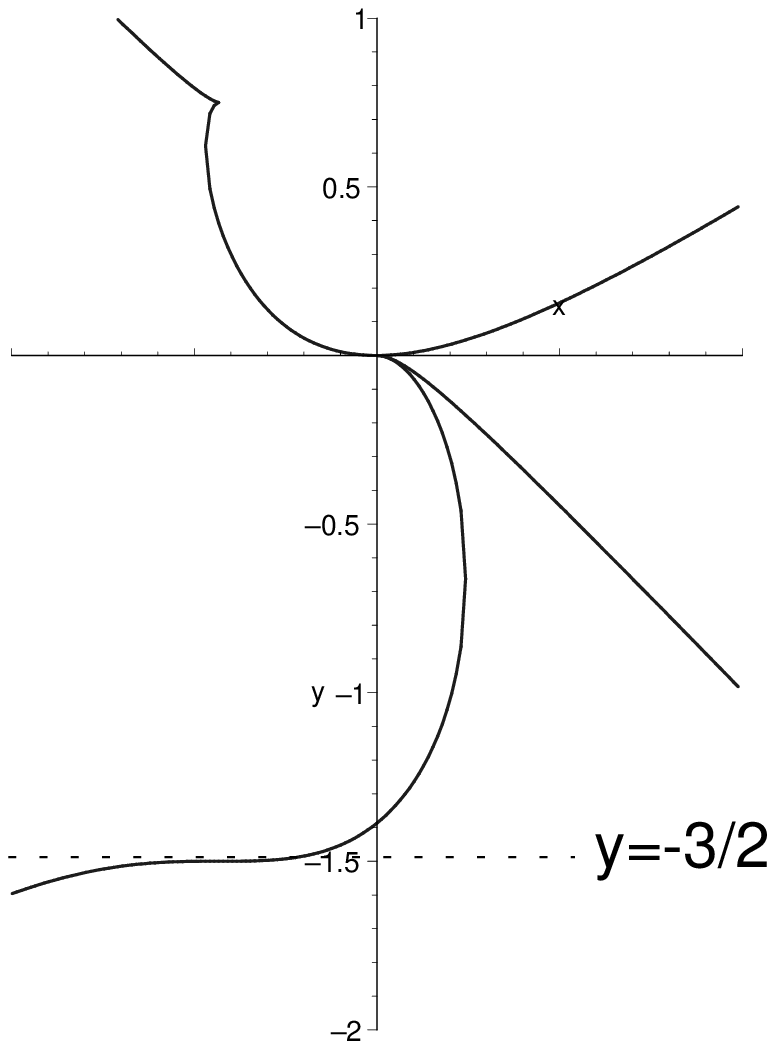}}
\put(20,130){\includegraphics{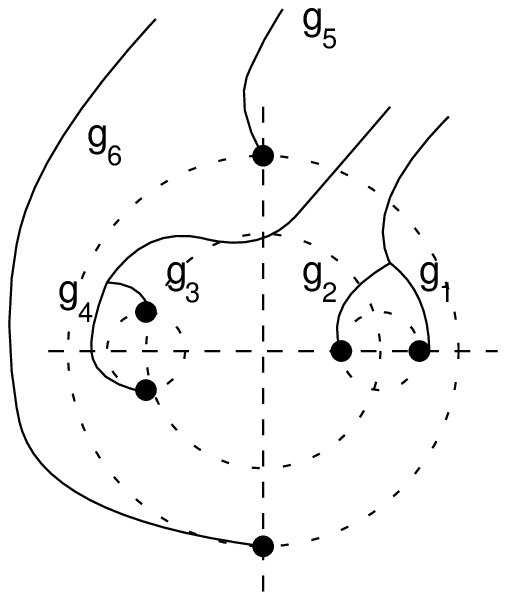}}
\end{picture}
\caption{Graph of $C_{\text{ex}}$ and its generators at $y=-\eps$}
\label{C(3,9)3A2}
\end{figure}

\noindent
To see the monodromy relation at the origin, we look at the  Puiseux
expansion of $(C_{\text{ex}},O)\cong C_{3,9}$. It has two components, a smooth
component $L$ and a
component $K$ of the (2,7)-cusp,
and they  are 
 parametrized as follows.
\begin{eqnarray*}
&L:  &x=\sqrt{4}\tau+\nnt, \quad y=\tau^2\\ 
&K:  &x=-\frac{i\sqrt{2}}{2}t^2-\frac{4\sqrt[4]{-6}}{9}
t^5+\nnt,\quad y=t^4
\end{eqnarray*}
\noindent
We take  generators $g_1,g_2,\dots,g_6$ of the fundamental group 
$\pi_1(L_{\eta_0}-C_{\text{ex}})$  as in Figure \ref{C(3,9)3A2}
where $\eta_0=-\eps$ and
$0<\varepsilon\ll 1$.
$g_5,g_6$ correspond to the points of $L$ and $g_1,\dots,g_4$
correspond to the points of $K$.
%\newline
Put $\om_1:=g_2g_1$, $\om_2:=g_4g_3$.
The monodromy relations at $O$ are given by:
\begin{eqnarray} \label{C3,9&3A2-O-L}
g_5=g_6,\quad g_6=(\om_2\om_1)g_5(\om_2\om_1)\inv \quad \text{(relation for L)}
\end{eqnarray} 
\begin{eqnarray} \label{C3,9&3A2-O-K1}
g_1=(g_5\om_2)g_3(g_5\om_2)^{-1},\quad g_2=(g_5\om_2)g_4(g_5\om_2)^{-1}
\end{eqnarray}
\begin{eqnarray} \label{C3,9&3A2-O-K2}
g_3=(g_5\om_2\om_1)g_2(g_5\om_2\om_1)^{-1},\quad 
g_4= (g_5\om_2\om_1^2)g_1(g_5\om_2\om_1^2)^{-1} 
\end{eqnarray}
%Put  $h:=g_5=g_6$. 
By (\ref{C3,9&3A2-O-K1}), we  reduce the generators to
$g_5,g_3,g_4$. 
Taking product of (\ref{C3,9&3A2-O-K1}), we get
% (\ref{C3,9&3A2-O-K1}) is equivalent to
\begin{eqnarray} \label{C3,9&3A2-K3}
\om_1=g_5\om_2 g_5\inv
%, \quad g_1=(g_5\om_2)g_3(g_5\om_2)^{-1}
\end{eqnarray}
The second relation of   (\ref{C3,9&3A2-O-L}) is equivalent to
\begin{eqnarray} \label{C3,9&3A2-O-L-b}
(g_5\om_2)^2=(\om_2 g_5)^2,\quad\text{or}\quad
(g_5g_4g_3)^2=(g_4g_3g_5)^2
\end{eqnarray} 
From  (\ref{C3,9&3A2-O-K1}) and (\ref{C3,9&3A2-K3}),
(\ref{C3,9&3A2-O-K2}) is rewritten as
\begin{eqnarray} \label{C3,9&3A2-O-K2b}
g_3=(g_5\om_2)^2\om_2 g_4\om_2\inv (g_5\om_2)^{-2},\quad 
g_4= (g_5\om_2)^2 \om_2^2 g_3\om_2^{-2} (g_5\om_2)^{-2} 
\end{eqnarray}
To read the monodromy relation at $y=-3/2$, which is a flex tangent
relation (with respect to $g_5,g_6,g_2$ in Figure \ref{C(3,9)3A2}).
\begin{eqnarray} \label{C3,9&3A2-flex}
g_3g_2g_3\inv=g_5
\end{eqnarray}
The singular pencil line $y=-9/4$ passes through two $A_2$ singularities,
thus the monodromy relations are braid relations, they are given by
\begin{eqnarray} \label{C3,9&3A2-2A2}
\{g_5,g_3\}=e,\quad\{g_5,g_4\}=e
\end{eqnarray} 
Under these braid relations, the relation (\ref{C3,9&3A2-flex}) reduces to 
the braid relation
\begin{eqnarray} \label{C3,9&3A2-A2-g3g4}
\{g_3,g_4\}=e
\end{eqnarray}
\begin{figure}[htb]
\setlength{\unitlength}{1bp}
\begin{picture}(255,130)(-100,0)
\put(-100,130){\includegraphics{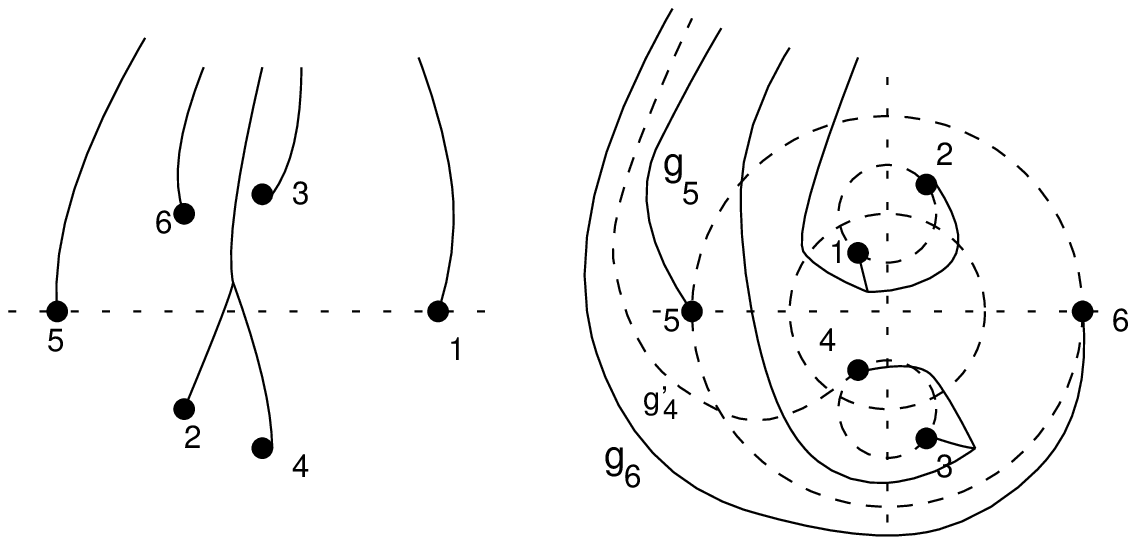}}
\end{picture}
\caption{Generators position at $y=(-9/4)^{+}$ (left) and at $y=0^{+}$ (right)}
\label{C3,9&3A2-halfturn}
\end{figure}
\noindent
To see the monodromy relation at $y=3/4$ first we show the position of
the generators after a half turn, see the right hand side of  Figure 
\ref{C3,9&3A2-halfturn}. The local equation of the  $A_2$ singularity at 
this level is $x^3+y^2=0$. Thus the monodromy relations are given by
\begin{eqnarray} \label{C3,9&3A2-A2}
g_1=g_4',\,\,\,\{g_5,g_1\}=e\quad
\text{where}\, g_4':=g_5(g_3\inv g_4 g_3)g_5\inv
\end{eqnarray} 
\noindent
See Figure 
\ref{C3,9&3A2-halfturn}.
\noindent
The first relation of (\ref{C3,9&3A2-A2}) follows from the braid relations.
The relations (\ref{C3,9&3A2-O-K2b}) and the second relation of 
(\ref{C3,9&3A2-A2}) reduces to the commuting relations:
$[(g_5g_4g_3)^2,g_4]=[(g_5g_4g_3)^2,g_3]=e$ which is equivalent to
\begin{eqnarray}
(g_5g_4g_3)^2=(g_4g_3g_5)^2=(g_3g_5g_4)^2
\end{eqnarray}
Thus the affine fundamental group $\pi_1(\bfC^2-C_{\text{ex}})$ is
isomorphic to
\begin{eqnarray}\label{affine-fund}
\langle g_3,g_4,g_5 |
\{g_5,g_3\}=\{g_5,g_4\}=\{g_3,g_4\}=e,
(g_5g_4g_3)^2=(g_4g_3g_5)^2=(g_3g_5g_4)^2\rangle
\end{eqnarray} 
Since $C_{\text{ex}}$ intersects transversely with the line at infinity
$L_{\infty}$, the generic  Alexander polynomial $\De(t)$ can be computed by Fox
calculus
\cite{Fox}. Let $p:\tilde X\to \bfP^2-(C_{\text{ex}}\cup L_\infty)$ be the infinite
cyclic covering corresponding to the kernel of the Hurewicz
homomorphism 
\[
\xi: \pi_1(\bfP^2-(C_{\text{ex}}\cup L_\infty))\to H_1(\bfP^2-(C_{\text{ex}}\cup
L_\infty);\bfZ)\cong \bfZ\]
 and put 
$\La:=\bfQ[t,t\inv]$.
Then the knot polynomials are given by 
  $\Delta(t)=\De_1(t)=(t^2-t+1)^2$ and $\De_2(t)=t^2-t+1$. 
Thus as $\La$-module, we have
\[H_1(\tilde X,\bfQ)\cong \La/(t^2-t+1)\oplus \La/(t^2-t+1)\]
Thus this implies $H_1(M;\bfQ)\cong \bfQ^4$. We can also show that 
$H_1(M;\bfZ)\cong \bfZ^4$ by computing the commutator subgroup
$D(\pi_1(\bfC^2-C_{\text{ex}}))$ using Reidmeister-Schreier method (\cite{M-K-S}).
%Thus we can claim $\pi_1(\bfP^2-C_{\text{ex}})\not\cong\Gr$.

To obtain the projective fundamental group, we add to
(\ref{affine-fund}) the big circle relation, which is given by
\begin{eqnarray} \label{C3,9&3A2-K4}
(g_5g_4g_3)^2=e
\end{eqnarray}
Relations $(g_4g_3g_5)^2=(g_3g_5g_4)^2=e$ follows from 
(\ref{C3,9&3A2-K4}). Thus we get
\[
\pi_1(\bfP^2-C_{\text{ex}}) \cong \langle g_3,g_4,g_5 |
\{g_5,g_3\}=\{g_5,g_4\}=\{g_3,g_4\}=e,(g_5g_4g_3)^2=e\rangle
\]
This completes the proof of Theorem \ref{main-theorem}.
We remark here that $\pi_1(\bfP^2-C_{\text{ex}})$ is much bigger than
$\bfZ_2*\bfZ_3$
by Corollary \ref{bigger} and the canonical surjection is 
given by 
identifying $g_3=g_5$.

\section{Non-tame torus curves}

\subsection{Braid group $B_4(\bfP^1)$}
First we recall that the braid group of $n$ strings in $\bfP^1$,
which is denoted by $B_n(\bfP^1)$, has the usual generators
$g_1,\dots, g_{n-1}$ and it has the representation (see for example
\cite{FB}):
\begin{eqnarray}\label{braid}
&g_ig_j=g_jg_i\quad &|i-j|\ge 2\\
&g_ig_{i+1}g_i=g_{i+1}g_ig_{i+1},\quad
&1\le i\le n-2\\
&g_1g_2\cdots g_{n-2}g_{n-1}^2g_{n-2}\cdots g_1=e
\end{eqnarray}
In particular, $B_4(\bfP^1)$ is generated by three elements
$g_1,g_2,g_3$ with relations:
\begin{eqnarray}\label{graid4}
&g_1g_3=g_3g_1,\quad & \{g_1,g_{2}\}=\{g_2,g_{3}\}=e\\
&g_1g_2g_3^2g_2g_1=e
\end{eqnarray}
Recall that $\bfZ_2*\bfZ_3$ is generated by two elements $\rho,\xi$
which satisfies the relations
$\{\rho,\xi\}=e, (\rho\xi)^3=e$.
The correspondence $g_1,g_3\mapsto \xi$ and $g_2\mapsto \rho$
defines a surjective homomorphism
$\Psi: B_4(\bfP^1)\to \bfZ_2*\bfZ_3$.
It is known (and easy to show) that $\Ker\Psi$ is not trivial.

%\vspace{.2cm}
\subsection{\bf  Example}
We are ready to  give an example which show that
the property $\pi_1(\bfP^2-C)\not \cong \bfZ_2*\bfZ_3$
is not so exceptional  for non-tame torus curves.
As an example,we consider a tame torus sextics with three $E_6$
singularities.  Note that such a curve is elliptic.
This curve can be degenerated into  rational curves in two ways. First
degeneration is  to combine two $E_6$ to make $B_{3,8}$ so that the
configuration  of singularities is $\{B_{3,8},E_6\}$ (\cite{Pho}). This
degeneration can be done in  the tame torus curves. As we have seen in (III),
the fundamental groups is unchanged by this degeneration. Another
degeneration is to put one
$A_1$
singularity outside of the conic $C_2$.
In \cite{Oka-sextics}, it is shown that the moduli of the 
sextics of torus type
with three $E_6$ is one-dimensional and it can be parametrized as 
\[ C_s: f(x,y,s)=(y^2+x^2-2x)^3+s(y^2-x^2)^2(x-1)^2=0\]
for $s\in \bfC^*$
 and 
$\pi_1(\bfP^2-C_s)\cong \bfZ_2*\bfZ_3$ for $s\ne 27, 0$
and $C_{27}$ is the unique sextics which obtain an $A_1$ singularity. 
The
 fundamental group $\pi_1(\bfP^2-C_{27})$
changes by this degeneration.
In fact, we have
\begin{Theorem}\label{non-tame}
$\pi_1(\bfP^2-C_{27})\cong B_4(\bfP^1)$ but the generic Alexander polynomial
remain 
unchanged.
\end{Theorem}

\begin{proof}
 It is  H. Tokunaga who has first observed that 
$\pi_1(\bfP^2-C_{27})\not\cong \pi_1(\bfP^2-C_s)\cong\bfZ_2*\bfZ_3$ with $s\ne 27$. His 
method is to use the existence of a certain finite covering \cite{Tokunaga}.

\begin{figure}[htb]
\setlength{\unitlength}{1bp}
\begin{picture}(255,160)(-80,0)
\put(-160,160){\includegraphics{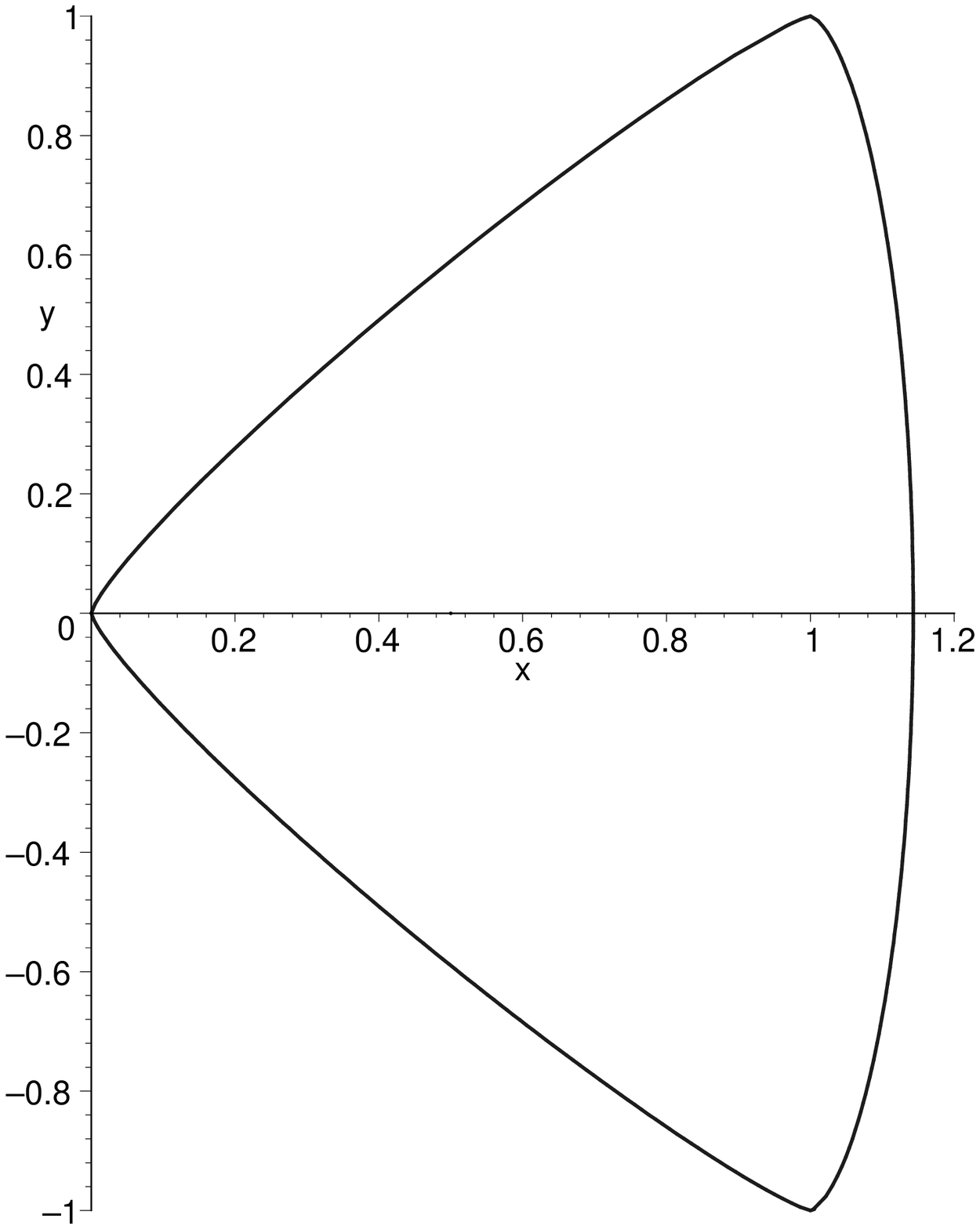}}
\put(-30,160){\includegraphics{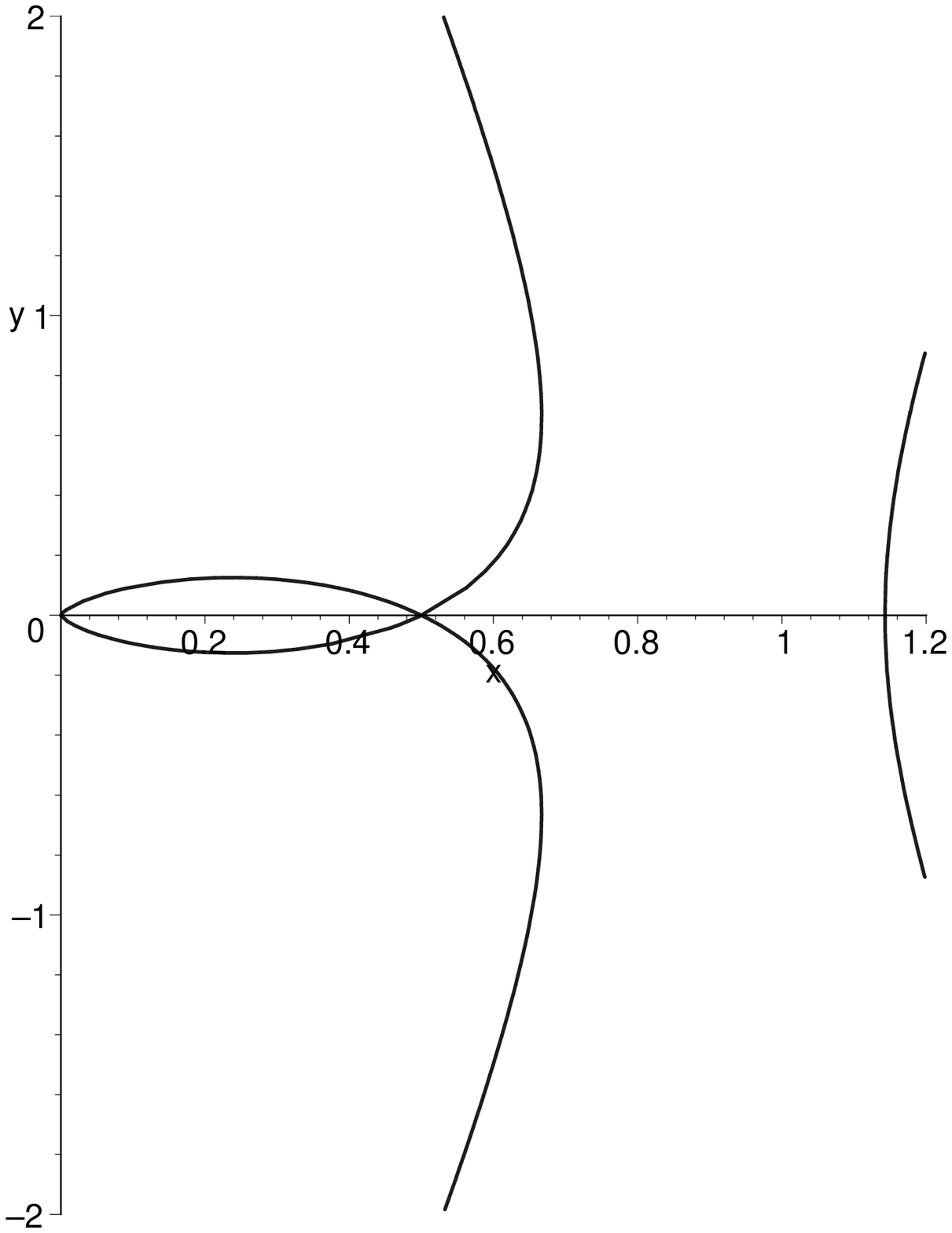}}
\put(80,160){\includegraphics{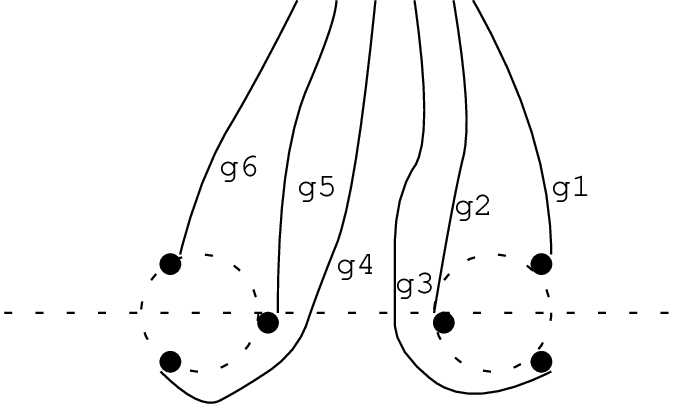}}
\end{picture}
\vspace{-2cm}
\caption{Graphs of $f(x,y)=0$, $f(x,yi)=0$ and generators at $x=1-\eps$}
\label{Graph-Quintic}
\end{figure}
\noindent
For the proof, we use the pencil $x=\eta, \eta\in \bfC$.
First
the discriminant polynomial of $\De_yf$ is given by
\begin{eqnarray*}
\De_y f(x,s)&=&c x^9 ((1+s)x^3-(6+2s) x^2+(12+s) x-8) (x-1)^{16} s^4
 ((27+2s) x-2 s)^2\\
\De_y f (x,27)&=& c x^9 (7x-8)(2x-1)^2(3x-2)^2(x-1)^{16}
\end{eqnarray*}
Observe that $(1+s)x^3-(6+2s) x^2+(12+s) x-8=0$ has three real roots 
for $s\ge 27$, where two of them make the multiple root $x=1/2$ for $s=27$.
Thus we have 5 (respectively 6) singular pencil lines for $s=27$
(resp. for $s>27$):
$x=0, 1/2, 2/3,1, 8/7$. In $C_{s}, s>27$, the node disappears and the
 singular line $x=1/2$ 
 splits into
two simple tangent lines $x=\alpha,\beta$ with $\alpha<1/2<\beta$. 
\newline
Hereafter
we consider the case $s=27$.
%The pencil $y=\eta$ is not suitable in our 
%computation as it has two non-real singular pencil. 
The line 
$x=1/2$ is passing the node $(1/2,0)\in C_{27}$, $x=2/3$ is a bitangent line
and $x=8/7$ 
is a simple tangent line.
The graphs of $f(x,y,27)=0$ and $f(x,y i,27)=0$
provide us the necessary informations.
We take 6 generators $g_1,g_2,g_3, g_4,g_5,g_6$ on the pencil line
$x=1-\eps$, $0<\eps\ll 1$ as in Figure \ref{Graph-Quintic}.

Put $\omega_1:=g_3g_2g_1$ and $ \omega_2:=g_6g_5g_4$. The big circle 
relation
is $ \omega_2 \omega_1=e$.
The monodromy relation at $x=1$  is given by 
\begin{eqnarray}
&g_1=\omega_1 g_2\omega_1\inv&,\quad g_2=\omega_1 g_3\omega_1\inv
\label{(1,1)}\\
&g_4=\omega_2g_5\omega\inv&,\quad g_5=\omega_2 g_6 \omega_2\inv
\label{(1,-1)}
\end{eqnarray}
At $x=8/7$, we get a simple tangent relation:
\begin{eqnarray}\label{what}
\omega_1\inv g_3\omega_1=\omega_2\inv g_6\omega_2
\end{eqnarray}
The line  $x=2/3$ is a bi-tangent line  and  the relation is:
\begin{eqnarray} \label{two-simple}
g_1=g_6,\quad g_3=g_4
\end{eqnarray}
Thus we can eliminate $g_4,g_6$ from generators.
At $x=1/2$, we get the commuting relation:
\begin{eqnarray} \label{commute}
g_1g_4=g_4g_1\quad\text{or},\quad g_1g_3=g_3g_1
\end{eqnarray}
Putting $\Omega:=g_5g_3g_2g_1$,
the monodromy relations at $x=0$ is
given by
\begin{eqnarray} 
&
(g_2g_1)\inv g_1 (g_2g_1)=g_5,&\label{0-1}\\
&g_1\inv g_2g_1=\Omega (g_2g_1)\inv g_1 (g_2g_1)\Omega\inv,& 
\quad
g_3=\Omega g_1\inv g_2g_1\Omega\inv\label{0-2}
 \end{eqnarray}
Now we can  rewrite (\ref{(1,1)}) using (\ref{commute})
as the braid relations:
%$g_1\omega=g_1g_3g_2g_1=g_3g_1g_2g_1$ and 
%$\omega_1g_2=g_3g_2g_1g_2$. This is equivalent to
%$g_2g_1g_2=g_1g_2g_1$.   Similarly  we get 
%$g_3g_2g_3=g_2g_3g_2$. Thus under (\ref{commute}), (\ref{(1,1)})
%and (\ref{(1,-1)}) are equivalent to 
\begin{eqnarray}\label{cusp}
g_1g_2g_1=g_2g_1g_2,\quad g_2g_3g_2=g_3g_2g_3
\end{eqnarray}
%The first  relation of  (\ref{(1,-1)})  reduces  to a trivial 
%one  and the second reduces to $g_5=g_2$ under
%$(\ref{cusp}) and  (\ref{two-simple}) and $\omega_2=g_6g_5g_4$ reduces to
The second relation of (\ref{0-1}) reduces to 
$g_5=g_2$.
Thus we can eliminate $g_5$  and we take $g_1,g_2,g_3$ as generators.
They satisfy (\ref{cusp}) and (\ref{commute}).
The relation
$\omega_2\omega_1=e$ reduces to 
\begin{eqnarray}\label{infty}
g_1g_2g_3^2g_2g_1=e
\end{eqnarray}
and we can see easily that other relations  follow  from
(\ref{cusp}),(\ref{commute}) and (\ref{infty}).
%\[g_3\omega_1^2\overset{\ref{commute}}{=}g_3g_3g_2g_3g_1g_2g_1
%\overset{\ref{cusp}}=
%g_3g_2g_3g_2g_1g_2g_1^2\overset{\ref{cusp}}=g_3g_2g_3g_1g_2g_1g_1
%\overset{\ref{commute}}=
%\omega_1^2g_1\]
 Thus we have proved 
\[
\pi_1(\bfP^2-C_{27})\cong \langle
g_1,g_2,g_3|\{g_1,g_2\}=\{g_2,g_3\}=e,g_1g_3=g_3g_1,(\ref{infty})\rangle
\cong B_4(\bfP^1)\]
The Alexander polynomial of $C_{27}$ is equal to $t^2-t+1$.
This can be shown by the exact same computation as in \cite{Oka-sextics}
or a direct Fox calculus \cite{Fox} from the above relations.
This completes the proof.
\end{proof}
Consider the degeneration $C_s,~s\to 27+0$.
Then we can see immediately that the monodromy relation at $x=1/2$
splits into   two simple tangent relations $x=\alpha$ and $x=\beta$ which
gives  the relation $g_3=g_1$ in the place of $g_1g_3=g_3g_1$.
This is the geometrical interpretation of the surjective homomorphism:
$\Psi:B_4(\bfP^1)\to \bfZ_2*\bfZ_3$.

\end{document}